\documentclass[10pt, journal]{IEEEtran}

\usepackage{epsfig,color,amsmath,cite}
\usepackage{amsthm} %defined already in ieeeconf
\usepackage{amsmath}    %For theorems
\usepackage[T1]{fontenc}
\usepackage[utf8]{inputenc}
\usepackage{bm}
\usepackage{epstopdf}
\usepackage{amssymb}
\usepackage{url}
\usepackage{enumitem} %defined already in ieeeconf
\usepackage{multirow}
\usepackage{hhline}
\usepackage{booktabs}
\usepackage{mathtools}
\usepackage{makecell}
\usepackage[linesnumbered,boxed,commentsnumbered,ruled,vlined,longend]{algorithm2e}
\usepackage{comment}

\DeclareMathOperator*{\minimize}{minimize}

\DeclareMathOperator*{\subjectto}{subject\ to}

\makeatother
\DeclareMathAlphabet\mathbfcal{OMS}{cmsy}{b}{n}

\newtheorem{mydef}{Definition}

\newtheorem{myrem}{Remark}

\newtheorem{myprs}{Proposition}



\usepackage{stackengine}

\newcommand{\mat}[1]{\boldsymbol{#1}}

\newcommand{\bmat}[1]{\begin{bmatrix} #1 \end{bmatrix}}

\providecommand{\mA}{\ensuremath{\mat{A}}}

\providecommand{\mC}{\ensuremath{\mat{C}}}

\providecommand{\mI}{\ensuremath{\mat{I}}}

\providecommand{\mL}{\ensuremath{\mat{L}}}

\providecommand{\mO}{\ensuremath{\mat{O}}}
\providecommand{\mP}{\ensuremath{\mat{P}}}

\providecommand{\mY}{\ensuremath{\mat{Y}}}
\providecommand{\mZ}{\ensuremath{\mat{Z}}}

\newcommand{\m}{\boldsymbol}
\allowdisplaybreaks[4]
\usepackage[colorlinks = true,
linkcolor = blue,
urlcolor  = blue,
citecolor = blue,
anchorcolor = blue]{hyperref}
\newcommand{\mbb}[1]{\mathbb{#1}}

\usepackage[framemethod=TikZ]{mdframed}
\mdfdefinestyle{MyFrame}{%
	linecolor=black,
	outerlinewidth=1.25pt,
	roundcorner=1.25pt,
	innerrightmargin=5pt,
	innerleftmargin=5pt,}
	
\usepackage[noabbrev]{cleveref}

\usepackage{mathtools}

\DeclarePairedDelimiter\abs{\lvert}{\rvert}%
\DeclarePairedDelimiter\norm{\lVert}{\rVert}%

\makeatletter
\let\oldabs\abs
\def\abs{\@ifstar{\oldabs}{\oldabs*}}
\let\oldnorm\norm
\def\norm{\@ifstar{\oldnorm}{\oldnorm*}}
\makeatother

\usepackage[english]{babel}
\usepackage[utf8]{inputenc}
\usepackage[super]{nth}

\usepackage{graphicx}
\usepackage{float}
\usepackage[caption = false]{subfig}

\usepackage{array}
\usepackage{threeparttable}

\usepackage[english]{babel}
\usepackage[utf8]{inputenc}
\usepackage[super]{nth}

\newcommand{\linf}{\mathcal{L}_{\infty}}

\def\newqed{{\null\nobreak\hfill\color{black}\ensuremath{\blacksquare}}}

\SetKwRepeat{Do}{do}{while}%

\title{Where Should Traffic Sensors Be Placed on Highways?}

\author{Sebastian A. Nugroh$\text{o}^{\ast}$, Suyash C. Vishno$\text{i}^{\dagger}$,  Ahmad F. Tah$\text{a}^{\ddagger,\diamond}$, Christian G. Claude$\text{l}^{\dagger}$, and Taposh Banerje$\text{e}^{\star}$ \vspace{-0.5cm}
	\thanks{
		
		$^{\ast}$Department of Electrical Engineering and Computer Science, University of Michigan, 1301 Beal Ave., Ann Arbor, MI 48109 (snugroho@umich.edu).

		$^\dagger$Department of Civil, Architectural, and Environmental 
		Engineering, The University of Texas at Austin, 301 E. Dean Keeton St. Stop C1700, Austin, TX 78712 (scvishnoi@utexas.edu,  christian.claudel@utexas.edu).
		
		$^\ddagger$Department of Civil and Environmental  Engineering, Vanderbilt University,  2201 West End Ave, Nashville, TN 37235 (ahmad.taha@vanderbilt.edu).		
		
		$^{\star}$Department of Electrical and Computer Engineering, The University of Texas at San Antonio, 1 UTSA Circle, San Antonio, TX 78249 (taposh.banerjee@utsa.edu).
		
		$^\diamond$Corresponding author. 
		
		This work is partially supported by the Valero Energy Corporation and National Science Foundation (NSF) under Grants 1636154, 1728629, 1739964, 1917164, 2152928, and 2152450.}
}

\begin{document}

\maketitle

\setlength{\abovedisplayskip}{3.5pt}
\setlength{\belowdisplayskip}{3.5pt}
\setlength{\abovedisplayshortskip}{3.2pt}
\setlength{\belowdisplayshortskip}{3.2pt}

\newdimen\origiwspc%
\newdimen\origiwstr%
\origiwspc=\fontdimen2\font% original inter word space
\origiwstr=\fontdimen3\font% original inter word stretch

\fontdimen2\font=0.64ex% inter word space
%\fontdimen3\font=0.66ex% inter word stretch

\begin{abstract}
This paper investigates the practical engineering problem of traffic sensors placement on stretched highways with ramps. Since it is virtually impossible to install bulky traffic sensors on each highway segment, it is crucial to find placements that result in optimized network-wide, traffic observability. Consequently, this results in accurate traffic density estimates on segments where sensors are \textit{not} installed. The substantial contribution of this paper is the utilization of control-theoretic observability analysis---jointly with integer programming---to determine traffic sensor locations based on the nonlinear dynamics and parameters of traffic networks. In particular, the celebrated asymmetric cell transmission model is used to guide the placement strategy jointly with observability analysis of nonlinear dynamic systems through Gramians. Thorough numerical case studies are presented to corroborate the proposed theoretical methods and various computational research questions are posed and addressed. The presented approach can also be extended to other models of traffic dynamics.
\end{abstract}


\vspace{-0.2cm}

\section{Motivation and Paper Contributions}\label{sec:Introduction}
%Intro + Literature Review + Paper Contributions + Organization + Notation (as in the ACC paper) \cite{nugroho2018journal}. 

\IEEEPARstart{W}{ith} the development of intelligent transportation systems technologies, numerous sensing methods for traffic data collection have become popular~\cite{Lovisari2016, Gentili2012}. Fixed sensors such as induction loops and magnetometers allow network operators to obtain high-quality measurements of vehicle density besides other information such as vehicle speeds and flow~\cite{Lovisari2016}. Sensors, however, are generally expensive to install and maintain which makes them infeasible for installation throughout the network and covering all segments. This poses a problem for traffic management operations and controls that require knowledge of the traffic state on all segments. Such operations include traffic control tasks such as ramp metering~\cite{Papageorgiou2002,Carlson2010a}, and variable speed limit control~\cite{Hegyi2010,Carlson2010a}. While it may only be feasible to collect data from a fixed number of segments in the network, it is still possible to obtain accurate estimates of the state of traffic on all the segments if the sensors are strategically placed. 

The problem of placing traffic sensors on highway networks has been divided into two main categories in the literature, one involving estimation and the other involving observability. Under the former category, the objective is to determine the optimal placement of sensors that minimizes the estimation error for unmeasured quantities such as travel time~\cite{Kianfar2009,Jovanovic2018,Gentili2018}, OD matrix~\cite{Mitsakis2017,Owais2019}, link flows\cite{Mehr2018}. Under observability, the literature is further divided based on the type of observability that is considered, full or partial observability. A fully observable system is one in which all the states are observable given the available measurements from the sensors. Different observability problems considered in this category are link flow observability~\cite{Hu2009, Ng2012, Bianco2001,Shao2016,Salari2019}, route flow observability~\cite{Castillo2008a,Rindali2017} and OD flow observailibty~\cite{Yang1998, Gan2005, Minguez2010}. The general idea is to determine a sensor placement configuration that ensures full observability of the system. Partial observability implies that not all states are observable given the set of available sensor measurements. Under this category, the literature focuses on determining such sensor locations that maximize the number of observable states while obeying some budget constraints. This concept is attractive in cases where the number of sensors required to achieve full observability is too large. Interested readers can refer to~\cite{Gentili2012, Viti2014, Castillo2015, Salari2019} and references therein for a detailed literature review of the aforementioned categories of sensor placement problem. In the current work, we focus on determining the optimal sensor placement to achieve full observability of the system while maintaining a balance between the number of sensors and the degree of observability of the system. Note that the degree of observability is a quantitative measure of the quality of state estimates that can be obtained by using a certain sensor placement configuration and is not related to the idea of partial observability. Unlike the aforementioned studies that utilize the relationship between the various flows in the network such as link and path flows, this study uses a traffic dynamics model to determine the relationship between various state variables. Also, the states considered here are traffic densities instead of flows.

A study that is closer to the work presented in this paper is~\cite{Contreras2016} which also considers a traffic dynamics model and applies a control theoretic approach to determine optimal sensor locations on {unconnected} highway segments. This approach is extended to study sensor placement on complex traffic networks in \cite{Agarwal2016}, and then adopted in \cite{Contreras2018} for studying the observability of highway traffic with Lagrangian sensors. These studies linearize the nonlinear traffic dynamics around a steady-state traffic flow. The major drawback of such approaches is that the linearized dynamics are valid only around the specific states. {Furthermore, these studies consider the Greenshield's fundamental diagram~\cite{Greenshields1935} which is inferior in modeling traffic flows when compared with the triangular and trapezoidal fundamental diagrams which are often considered in the implementation of the \textit{cell transmission model} (CTM)~\cite{Gomes2008, Daganzo1994}.} \cite{Lovisari2016} does deal with optimal sensor placement and density reconstruction while considering the CTM model but to simplify things it still linearizes the model and thus has the same drawback as above.

There exist some other papers studying the observability properties of traffic systems considering a traffic dynamics model, albeit do not formally study the sensor placement. The observability properties related to the switching modes model\cite{Munoz2003}, which is a piecewise-linearized version of the CTM\cite{Daganzo1994} is studied in \cite{Munoz2006}. In \cite{Ferrara2016}, the authors propose a ramp-metering controller based on a switching modes scheme and study the observability properties of the considered freeway system. A new piecewise affine system model based on CTM is developed in \cite{Guo2018}, where its observability is investigated. These studies however, linearize the traffic dynamics and suffer the same drawback.

This paper investigates the sensor placement problem for highway networks with ramps having nonlinear traffic dynamics. The traffic model is built using the \textit{asymmetric CTM} (ACTM)~\cite{Gomes2004, Gomes2006, Gomes2008} (a close variant of the CTM) and the triangular fundamental diagram, and thus is also an improvement over \cite{Contreras2016, Agarwal2016, Contreras2018} in terms of modeling traffic dynamics, to remove the dependence of optimal sensor placements on assumptions about traffic states. The traffic sensor placement is addressed via observability analysis based on the nonlinear traffic model, which is different from the ones used in the aforementioned papers. 

%Unlike linear systems, determining observability for nonlinear system, even without incorporating sensor selection, is not as straightforward. 
There exist several approaches in the literature to quantify observability for nonlinear systems. The two most prevalent approaches include utilizing the concept from differential embedding and Lie derivative \cite{Letellier2005,Whalen2015}, while the other is based on empirical observability Gramian \cite{LALL19992598}. In the context of sensor placement for nonlinear systems, four different approaches have been proposed recently. The first approach leverages empirical observability Gramian to find the best set of sensors that maximizes some metrics on the Gramian matrix---see \cite{SERPAS2013105,qi2015optimal,singh2006sensor}. Another approach is proposed in \cite{Haber2017} where the observability Gramian around a certain initial state is constructed through a moving horizon estimation (MHE) framework. As such, the sensor selection problem is posed as an integer programming (IP) problem which objective is maximizing the logarithmic determinant of the Gramian. The third approach, developed in \cite{Bopardikar2019}, introduces a randomized algorithm for dealing with the sensor placement problem and accordingly, theoretical bounds for eigenvalue and condition number of observability Gramian are proposed. The last and most recent approach is established in \cite{Nugroho2019CDC} where the authors make use of the numerous observer designs for some classes of nonlinear systems posed as semidefinite programs (SDP). The resulting problem, which combines the sensor selection \textit{and} observer gain synthesis, is posed as a nonconvex mixed-integer SDP, which can be cast into a convex one through some relaxation/reformulation techniques.

%\comment{Suyash: What about that other paper which considers switched system model and then finds placements? Why didn't it get a mention here?} 
%
%\comment{Suyash: Second, you need to articulate how there's a research gap in the sense that the literature (Agrawal and that other switched system paper) does not consider what we consider here or something like that. After that, you can explain the paper's contribution. Please read the abstract I wrote and the title, and then tie this up nicely with the paper's contributions. Finally, you need to make the connection between state estimation and sensor placement. In short, the next paragraphs need to be re-written. The above paragraphs are mostly fine.}

Herein, based on the nonlinear traffic model developed using ACTM, we formulate the sensor placement problem using traffic observability analysis.
%construct the sensor placement problem through the MHE framework developed in \cite{Haber2017}. 
The resulting problem---categorized as convex IP---is solved via an integer branch-and-bound (BnB) algorithm and therefore, optimal traffic sensor locations can be obtained. Notice that the use of nonlinear traffic model and observability analysis result in traffic sensor locations that are valid for various traffic conditions.
The novelties made in this paper are listed as follows:  
\vspace{-0.05cm}
\begin{itemize}[leftmargin=*]
	\item We present a discrete-time nonlinear state-space model for stretched highway networks having multiple on- and off-ramps based on ACTM. In addition, we analytically compute the corresponding Lipschitz constant for the nonlinear counterparts. The formulated state-space model along with Lipschitz constant pave the way for implementing various control-theoretic approaches to solve control and estimation problems prevailing in highway traffic networks.
	\item We leverage the concept of observability for nonlinear systems to construct the traffic sensor placement problem, which is equivalent to maximizing the determinant and trace of observability Gramian matrix given the number of allocated traffic sensors. To the best of our knowledge, this is the first attempt to solve traffic sensor placement problem based on the degree of observability for the nonlinear traffic dynamics.
	\item We propose a corresponding robust observer framework for Lipschitz nonlinear discrete-time systems developed using the concept of $\mathcal{L}_{\infty}$ stability for traffic density estimation purpose. By using the optimal configuration of traffic sensors, the computation of stabilizing observer gain matrix can be cast as a SDP and as such, can be solved efficiently.
	\item We verify the effectiveness of our approach 
	to solve the traffic sensor placement problem through various case studies. Specifically, we \textit{(a)} study the performance and computational efficiency of determinant and trace functions in estimating the system's initial state, \textit{(b)} assess their impact on traffic density estimation performance under noisy conditions, and \textit{(c)} compare their results relative to randomized and uniform sensor placement strategies.
\end{itemize}  
\vspace{-0.05cm}

The remainder of the paper is organized as follows. In Section \ref{sec:nonlinearmodel}, we provide the mathematical modeling of stretched highway with ramps using ACTM, which result in a nonlinear state-space form. Next, Section \ref{sec:traffic_ssp} discusses our strategy for addressing traffic sensor placement based on system's observability. Section \ref{sec:observer_design} develops a simple robust observer design for traffic density estimation purpose. The proposed approach is extensively tested in Section \ref{sec:cases_tudy} for solving the traffic sensor placement problem and density estimation via numerous case studies. Finally, Section \ref{sec:summary} concludes the paper.   

\noindent {\textbf{Paper's Notation:}} \; Let $\mathbb{N}$, $\mathbb{R}$, $\mathbb{R}_{++}$, $\mathbb{R}^n$, and $\mathbb{R}^{p\times q}$ denote the set of natural numbers, real numbers, positive real numbers, and real-valued row vectors with size of $n$, and $p$-by-$q$ real matrices respectively. $\mathbb{S}^{m}$ denotes the set of symmetric matrices with the dimension of $m$. Specifically, $\mathbb{S}^{m}_{+}$ and $\mathbb{S}^{m}_{++}$ denotes the set of positive semi-definite and positive definite matrices respectively. For any vector $z \in \mathbb{R}^{n}$, $\Vert z\Vert_2$ denotes its Euclidean norm, i.e. 
$\Vert z\Vert_2 = \sqrt{z^{\top}z} $, where $z^{\top}$ is the transpose of $z$. 
%For two vectors $u,v \in \mathbb{R}^{n}$, then $\langle u,v\rangle = u^{\top}v$, which denotes the inner product between those two vectors. 
The symbol $\otimes$ denotes the Kronecker product where $\mathrm{det}(\m A)$ and $\mathrm{trace}(\m A)$ return the determinant and trace of matrix $\m A$. 
For simplicity, the notation '$*$' denotes terms induced by symmetry in symmetric block matrices. 
Tab. \ref{tab:notation} provides the nomenclature utilized in this paper.

%\noindent \textbf{Note about Supplemental Document:} The appendices of this paper, containing detailed mathematical derivations and proofs, are all included in a supplemental document attached with this manuscript submission. 

\setlength{\textfloatsep}{10pt}
\begin{table}[t!]
	\footnotesize	\renewcommand{\arraystretch}{1.3}
	\caption{Paper nomenclature: parameter, variable, and set definitions.}
	\label{tab:notation}
	\centering
	\begin{tabular}{||l|l||}
		\hline
		\textbf{Notation} & \textbf{Description}\\
		\hline
		\hline
		\hspace{-0.1cm}$\Omega$ & \hspace{-0.1cm}the set of highway segments on the stretched highway \\
		\hspace{-0.1cm} & \hspace{-0.1cm}$\Omega = \{ 1,2,\hdots,N \}$ , $N := \abs{\Omega}$ \\
		\hline
		\hspace{-0.1cm}$\Omega_I$ & \hspace{-0.1cm}the set of highway segments with on-ramps \\
		\hspace{-0.1cm}	&  \hspace{-0.1cm}$\Omega_I = \{ 1,2,\hdots,N_I \}$ , $N_I := \abs{\Omega_I}$ \\
		\hline
		\hspace{-0.1cm}$\Omega_O$ & \hspace{-0.1cm}the set of highway segments with off-ramps \\
		\hspace{-0.1cm}& \hspace{-0.1cm}$\Omega_O = \{ 1,2,\hdots,N_O \}$, $N_O := \abs{\Omega_O}$ \\ 
		\hline
		\hspace{-0.1cm}$\hat{\Omega}$ & \hspace{-0.1cm}the set of on-ramps, $\hat{\Omega} = \{ 1,2,\hdots,N_I \}$ , $N_I = |\hat{\Omega}|$\\
		\hline
		\hspace{-0.1cm}$\check{\Omega}$ & \hspace{-0.1cm}the set of off-ramps, $\check{\Omega} = \{ 1,2,\hdots,N_O \}$ , $N_O = |\check{\Omega}| $\\
		\hline
		\hspace{-0.1cm}$T$ & \hspace{-0.1cm}duration of each time step\\
		\hline		
		\hspace{-0.1cm}$l$ & \hspace{-0.1cm}length of each segment, on-ramp, and off-ramp\\
		\hline	
		
		\hspace{-0.1cm}$\rho_i [k]$ & \hspace{-0.1cm}traffic density on Segment $i \in \Omega$ at time $kT$, $k \in \mathbb{N}$ \\
		\hline
    	\hspace{-0.1cm}$q_i[k]$ & \hspace{-0.1cm}traffic flow from Segment $i \in \Omega $ into the next segment\\ 
		\hline
    	\hspace{-0.1cm}$\delta_i[k], \sigma_i[k]$ & \hspace{-0.1cm}demand and supply functions for Segment $i \in \Omega$\\
		\hline
		
		\hspace{-0.1cm}$\hat{\rho}_i[k]$ & \hspace{-0.1cm}traffic density on the on-ramp of Segment $i\in \Omega_I $ \\
		\hline
    	\hspace{-0.1cm}$r_i[k]$ & \hspace{-0.1cm}traffic flow into Segment $i \in \Omega_I $ from the on-ramp\\
		\hline
    	\hspace{-0.1cm}$\hat{r}_i[k]$ & \hspace{-0.1cm}traffic flow into the on-ramp of Segment $i \in \Omega_I $\\
		\hline
    	\hspace{-0.1cm}$\hat{\delta}_i[k], \hat{\sigma}_i[k]$ & \hspace{-0.1cm}demand and supply functions for the on-ramp of\\
    	\hspace{-0.1cm}	&  \hspace{-0.1cm}Segment $i \in \Omega_I$\\
		\hline		
		
		\hspace{-0.1cm}$\check{\rho}_i[k]$ & \hspace{-0.1cm}traffic density on the off-ramp of Segment $i\in \Omega_O$ \\
		\hline
    	\hspace{-0.1cm}$s_i[k]$ & \hspace{-0.1cm}traffic flow from Segment $i \in \Omega_O $ into the off-ramp\\
		\hline
    	\hspace{-0.1cm}$\check{\delta}_i[k], \check{\sigma}_i[k]$ & \hspace{-0.1cm}demand and supply functions for the off-ramp of\\
    	\hspace{-0.1cm}	&  \hspace{-0.1cm}Segment $i \in \Omega_O$\\
		\hline

    	\hspace{-0.1cm}$\check{s}_i[k]$ & \hspace{-0.1cm}traffic flow from the off-ramp of Segment $i \in \Omega_O $\\
		\hline
		
    	\hspace{-0.1cm}$f_{in}[k]$ & \hspace{-0.1cm}traffic wanting to enter Segment 1 of the highway\\
		\hline
    	\hspace{-0.1cm}$f_{out}[k]$ & \hspace{-0.1cm}traffic that can leave Segment $N$ of the highway\\
		\hline
    	\hspace{-0.1cm}$\hat{f}_i[k]$ & \hspace{-0.1cm}traffic wanting to enter the on-ramp of Segment $i \in \Omega_I $\\
		\hline
    	\hspace{-0.1cm}$\check{f}_i[k]$ & \hspace{-0.1cm}traffic that can leave the off-ramp of Segment $i \in \Omega_O $\\
		\hline

    	\hspace{-0.1cm}$\beta_i[k]$ & \hspace{-0.1cm}split ratio for the off-ramp of Segment $i \in \Omega_O$,\\
    	\hspace{-0.1cm}	&  \hspace{-0.1cm}where $\beta_i[k] \in [0,1]$\\
		\hline
    	\hspace{-0.1cm}$\xi_i[k]$ & \hspace{-0.1cm}occupancy parameter for the on-ramp of Segment $i \in \Omega_I$\\
    	\hspace{-0.1cm}	&  \hspace{-0.1cm}where $\xi_i[k] \in [0,w_c]$\\
    	\hline              
		\hspace{-0.1cm}$v_f$ & \hspace{-0.1cm}free-flow speed \\
		\hline
		\hspace{-0.1cm}$w_c$ & \hspace{-0.1cm}congestion wave speed \\
		\hline
		\hspace{-0.1cm}$\rho_m$ & \hspace{-0.1cm}maximum density  \\
		\hline
		\hspace{-0.1cm}$\rho_c$ & \hspace{-0.1cm}critical density \\
		\hline
	\end{tabular}
	\vspace{-0.0cm}
\end{table}
\setlength{\floatsep}{10pt}

\begin{figure}
	\centering
	\includegraphics[width=0.27\textwidth]{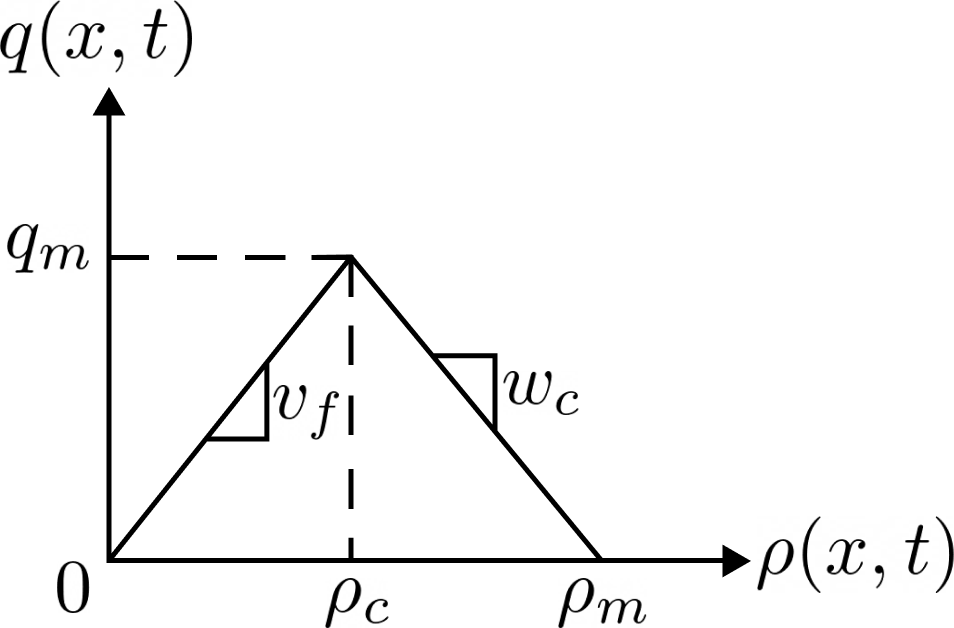}
	\vspace{-0.2cm}
	\caption{Triangular fundamental diagram.}
	\label{fig:Fund_diag}
	\vspace{-0.1cm}
\end{figure}

\section{Nonlinear Discrete-Time Modeling of Traffic Networks with ramps}\label{sec:nonlinearmodel}

This section presents the discrete-time modeling of traffic dynamics on a stretched highway with arbitrary number and location of ramps. 
To that end, here we utilize the \textit{Lighthill-Whitman-Richards} (LWR) Model \cite{Lighthill1955b,Richards1956} for traffic flow. 
%  which is expressed by a partial differential equation given as
% \begin{align}
%  \frac{\partial \rho (t,d)}{\partial t} + \frac{\partial q(t,d)}{\partial d} = 0 , \label{eq:LWRmodel}
% \end{align} 
In this paper, the relationship between traffic density and flux is given by the triangular-shaped fundamental diagram which has been extensively used in the literature \cite{Gomes2008}. It is depicted in Fig. \ref{fig:Fund_diag} and is constructed as
\begin{align}
\small\hspace{-0.2cm}q\left(\rho(t,d)\right) \hspace{-0.05cm}= \hspace{-0.05cm}
\begin{cases}\hspace{-0.05cm}
v_f\rho(t,d), &\;\mathrm{if} \; 0 \leq \rho(t,d) \leq \rho_c \\
\hspace{-0.05cm}w_c\left(\rho_m-\rho(t,d)\right),&\; \mathrm{if} \;\rho_c \leq \rho(t,d) \leq \rho_m.\\
\end{cases}\label{eq:triangular_model}
\end{align}
where $t$ and $d$ denote the time and distance; $\rho(t,d)$ denotes the traffic density (vehicles/distance) and $q(t,d)$ denotes the traffic flux (vehicles/time).
 To represent the traffic dynamics as a series of difference, state-space equations---a useful bookkeeping for the ensuing discussions---we discretize the LWR Model with respect to both space and time (this is also referred to as the Godunov discretization). This approach allows the highway of length $L$ to be divided into segments (cells) of equal length $l$ and the traffic networks model to be represented by discrete-time equations. These segments form both the highway and the attached ramps. Throughout the paper, the segments forming the highway are referred to as mainline segments. We assume that the highway is split into $N$ mainline segments.
 
To ensure computational stability, the Courant-Friedrichs-Lewy condition (CFL) given as ${v_f T}l^{-1}\leq 1$ has to be satisfied. Since each segment is of the same length $l$, then we have $\rho(t,d) = \rho(kT,l)$, where $k\in\mathbb{N}$ represents the discrete-time index. For simplicity of notation, from here on $\rho(kT,l)$ is simply written as $\rho[k]$.

As mentioned earlier, traffic is modeled using the ACTM which is originally given in~\cite{Gomes2004, Gomes2006, Gomes2008}. {The ACTM is a variant of the CTM that departs from the CTM in its treatment of asymmetric merge junctions such as the on-ramp-highway junctions. Unlike the CTM, it assumes separate allocations of the available space on the highway for traffic from each merging branch, which allows for comparatively simple flow conservation equations at those merges than the original CTM.} A variation to the original ACTM is introduced in the modeling of the ramps which here are treated as normal segments rather than point queues as in the original approach. The approach in this paper is similar to that in~\cite{Vishnoi2020}. Herein, we define new functions referred to as the demand function $\delta_i[\cdot]$ and the supply function $\sigma_i[\cdot]$ to simplify the ensuing expressions. The demand function equals the traffic flux leaving Segment $i$ through the highway assuming that the next segment has infinite storage. The supply function equals the traffic flux that can enter Segment $i$ through the highway assuming that the previous segment has infinite storage. These functions are constructed from the triangular FD given in \eqref{eq:triangular_model}.
The demand function $\delta_i[\cdot]$, with and without off-ramp, can be written as
%\begin{subequations}
	\begin{align}\label{eq:demand functions}
\hspace{-0.42cm}\small	\delta_{i}[k] \hspace{-0.05cm}= \hspace{-0.05cm}\begin{cases}
	\hspace{-0.05cm}\min \big(\bar{\beta_{i}}[k]v_{f}\rho_{i}[k],\bar{\beta_{i}}[k]v_{f}\rho_{c},\frac{\bar{\beta_{i}}[k]}{\beta_{i}[k]}\check{\sigma}_i[k]\big), \;\mathrm{if}\; i\in\Omega_O\\
	\hspace{-0.05cm}\min \big(v_{f}\rho_{i}[k],v_{f}\rho_c\big), \qquad\qquad\qquad\;\;\,\;\mathrm{if}\; i\in\Omega\setminus\Omega_O.
	\end{cases}
	\end{align}
%\end{subequations}
Here $\check{\sigma}_i[k]$ can be given as
\begin{align}
    \check{\sigma}_i[k]=\min\big(w_c(\rho_m - \check{\rho_i}[k]),v_f\rho_c\big).
\end{align} 
 In \eqref{eq:demand functions}, the split-ratio $\beta_i[k]$ relates the traffic flow from Segment $i$ into its off-ramp with the traffic flow from Segment $i$ into the next segment such that
%\begin{equation}
\begin{align}\label{eq:off-ramp relation}
    s_i[k]&=\beta_i[k](s_i[k]+q_i[k]) \implies
  s_i[k]=\frac{\beta_i[k]}{\bar{\beta_i}[k]}q_i[k].
\end{align}
%\end{equation}
where we define $\bar{\beta_i}[k]\overset{\Delta}{=}1-\beta_i[k]$ to simplify the equations. Since $\frac{\bar{\beta_{i}}[k]}{\beta_{i}[k]}v_f\rho_c>\bar{\beta}_iv_f\rho_c$, therefore in \eqref{eq:demand functions} we can also simply write $\check{\sigma}_i[k]=w_c(\rho_m - \check{\rho_i}[k])$.

Similarly, the supply function $\sigma_i[k]$, with and without on-ramp, can be represented as
%\begin{subequations}
    \begin{align}\label{eq:supply functions}
\small	\hspace{-0.43cm}\sigma_{i}[k] \hspace{-0.05cm}=\hspace{-0.05cm} \begin{cases} \hspace{-0.05cm}\min\big(w_{c}(\rho_{m}\hspace{-0.05cm}-\hspace{-0.05cm}\rho_{i}[k]),v_f\rho_c\big)-r_{i}[k],\,\qquad\hspace{0.2cm}\mathrm{if}\, i\in\Omega_I \\
 \hspace{-0.05cm}\min\big(w_{c}(\rho_{m}\hspace{-0.05cm}-\hspace{-0.05cm}\rho_{i}[k]),v_f\rho_c),\,\qquad\qquad\mathrm{if}\, i\in\Omega\setminus\Omega_I,
	\end{cases}
	\end{align}
%\end{subequations}
where $r_i[k]$ is described by the following equation
\begin{align}
    r_i[k]=\min\big(v_f\hat{\rho_i}[k],\xi_i(\rho_m-\rho_i[k]),\frac{\xi_i}{w_c}v_f\rho_c\big).\label{eq:on-ramp flow}
\end{align}
Here $v_f\hat{\rho_i}[k]$ is the demand of the on-ramp in free-flow, and $\frac{\xi_i[k]}{w_c}$ is the fraction of the flow that Segment $i$ can receive from its on-ramp.
The upstream and downstream flows for any mainline Segment $i$ are given as
\vspace{-0.0cm}
\begin{subequations}\label{eq:CTM_flow_characteristics}\begin{align}
q_{i-1}[k] &= \min \big( \delta_{i-1}[k],\sigma_i[k]\big) \label{eq:CTM_flow_characteristics_1} \\
q_{i}[k] &= \min\big(\delta_i[k],\sigma_{i+1}[k]\big), \label{eq:CTM_flow_characteristics_2}
\end{align}
\end{subequations}
%Since we use Godunov discretization with fixed time step $T$ and distance $l$, from \eqref{eq:speed-rho} and \eqref{eq:CTM_flow_characteristics} we obtain
%\begin{subequations}\label{eq:CTM_flow_characteristics_b}
%	\begin{align}
%	\phi_{i-1}[k] &= \frac{T}{l}q_{i-1}[k] = \frac{T}{l}\min\big(\delta_{i-1}[k],\,\sigma_i[k]\big) \label{eq:CTM_flow_characteristics_1b} \\
%	\phi_{i}[k] &= \frac{T}{l}q_{i}[k] = \frac{T}{l}\min\big(\delta_{i}[k],\,\sigma_{i+1}[k]\big). \label{eq:CTM_flow_characteristics_2b}
%	\end{align}
%\end{subequations}

% The flow conservation equations for the on-ramp and off-ramp can be written as follows
% \begin{subequations}\label{eq:Flow_conservation_ramps}
% \begin{align}
%     \hat{\rho_i}[k+1] &= \hat{\rho_i}[k] +\frac{T}{l}( \hat{r_i}[k] - r_i[k])\label{eq:Flow-conservation_on-ramp}\\
%     \check{\rho_i}[k+1] &= \check{\rho_i}[k] +\frac{T}{l}( s_i[k]-\check{s_i}[k]),\label{eq:Flow-conservation_off-ramp}
% \end{align}
% \end{subequations}
% where $\hat{r_i}[k]$ and $\check{s_i}[k]$ are given as
% \begin{subequations}\label{eq:ramp inflow outflow}
%     \begin{align}
%     \hat{r_i}[k]&=\min\big(w_c(\rho_m-\hat{\rho_i}[k]),v_f\rho_c,\hat{f_i}\big)\label{eq:on-ramp inflow}\\
%     \check{s_i}[k]&=\min\big(v_f\check{\rho_i}[k],v_f\rho_c,\check{f_i}\big),\label{eq:off-ramp outflow}
%     \end{align}
% \end{subequations}

The discrete-time flow conservation equation for a  mainline segment with both on- and off-ramp can then be written as %follows
% \begin{subequations}
% \label{eq:CTM_flow_conservation}
% \begin{align}
% \hspace{-0.2cm}\rho_i[k+1] &=\rho_i[k]+\frac{T}{l}\big(q_{i-1}[k]-q_{i}[k]\big)\\
% \hspace{-0.2cm}\hspace{-0.2cm}\rho_i[k+1] &=\rho_i[k]+\frac{T}{l}\big(q_{i-1}[k]+r_{i}[k]-q_{i}[k]\big)\\
% \hspace{-0.2cm}\rho_i[k+1] &=\rho_i[k]+\frac{T}{l}\big(q_{i-1}[k]-q_{i}[k]-s_{i}[k]\big)
% \end{align}
% \end{subequations}

%\begin{equation}
\begin{align}
\hspace{-0.2cm}\rho_i[k+1] &=\rho_i[k]+\frac{T}{l}\big(q_{i-1}[k]+r_{i}[k]-q_{i}[k]-s_{i}[k]\big)\label{eq:CTM_flow_conservation}
\end{align}
%\end{equation}

\begin{figure}[t]
	\centering
	\includegraphics[width=0.36\textwidth]{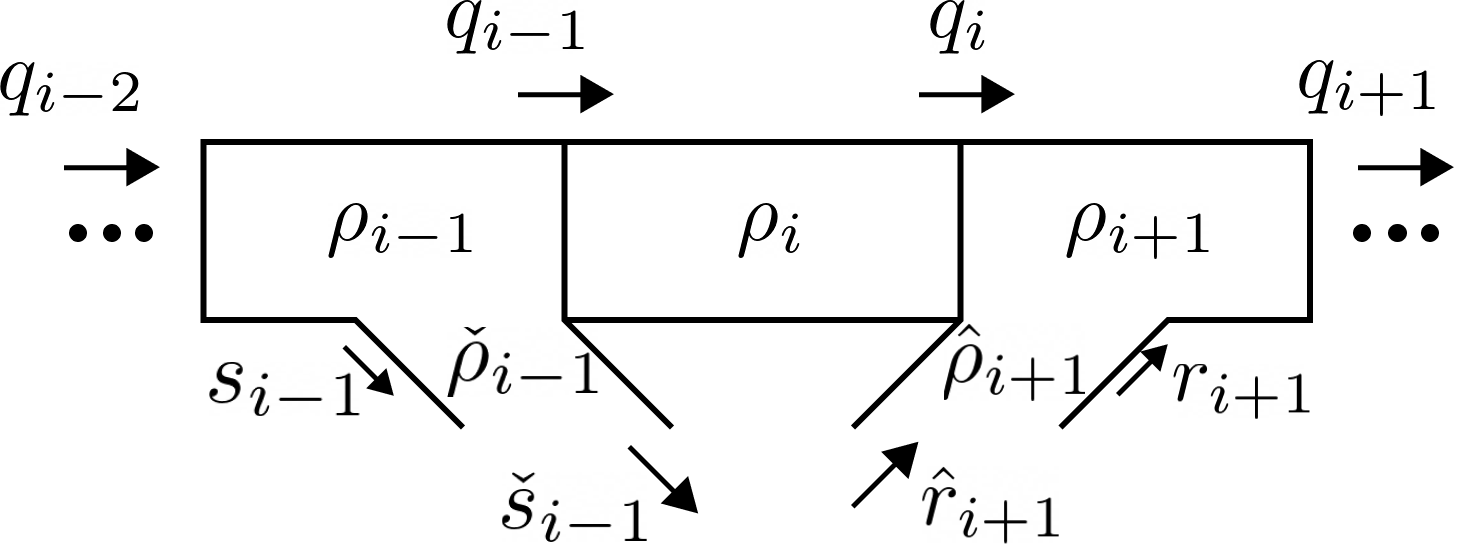}
	\vspace{-0.2cm}
	\caption{Consecutive segments of the highway with ramp connections.}
	\label{fig:Highway_Piece}
%	\vspace{-0.5cm}
\end{figure}

Equations for segments with only one or no ramp can be written by removing the respective flow terms. The flows and conservation equations for the ramps can also be written similar to the mainline segments. The detailed illustration for the above model is given in Fig. \ref{fig:Highway_Piece}.

\begin{figure}
	\centering
	\includegraphics[width=0.35\textwidth]{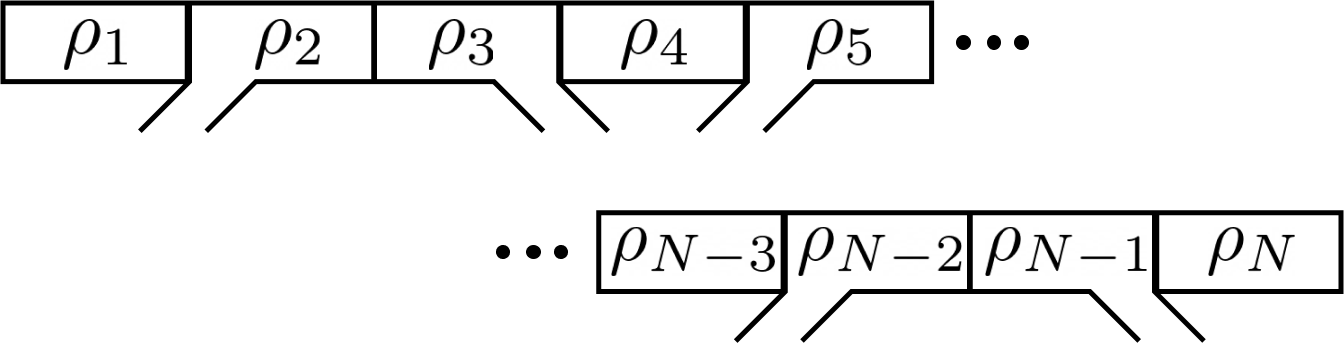}
	\vspace{-0.2cm}
	\caption{Schematic of the highway structure considered in this study.}
	\label{fig:Highway_Full}
	\vspace{-0.1cm}
\end{figure}

%\vspace{-0.24cm}
In the present study, we assume an arbitrary structure of the highway with respect to the positioning of the ramps on different mainline segments. In that, we assume that every mainline segment having no ramps is followed by a segment having an on-ramp, which is followed by a segment having an off-ramp, and which is followed by a segment having no ramps and this structure repeats throughout the highway. Additionally, we assume that the first and last mainline segments have no ramps. For example, a highway split into seven segments has an on-ramp on Segments 2 and 5, an off-ramp on Segments 3 and 6, and no ramps on Segments 1, 4 and 7. Fig. \ref{fig:Highway_Full} gives a schematic representation of the structure of the highway considered in this paper. Note that this structure of the highway is arbitrarily chosen for experimentation of the methodology presented in this work. It is easy to consider a different highway structure as per an actual site where the method may be used.

In this paper we assume that the upstream demand at Segment $1$ and downstream supply at Segment $N$ are known, denoted by $\delta_0[k]=f_{\mathrm{in}}[k]$ and $\sigma_{N+1}[k]=f_{\mathrm{out}}[k]$ respectively. Similarly, we assume that all the on-ramp demands $\hat{f}_i,i\in\Omega_I$ and off-ramp supplies $\check{f}_i,i\in\Omega_O$ are also known. These known demand and supply values are lumped into the input vector defined as $$\m u[k]=[f_{\mathrm{in}}\hspace{1mm} f_{\mathrm{out}}\hspace{1mm} \ldots \hspace{1mm}{\hat{f}_j}\hspace{1mm} \ldots \hspace{1mm}{\check{f}_l}\hspace{1mm} \ldots]^{\top}\in \mbb{R}^{2+N_I+N_O},$$ where $j\in\Omega_I$ and $l\in\Omega_O$.
The state vector can be defined as $$\m x[k] := [\rho_i[k]\hspace{1mm} \ldots\hspace{1mm}\hat{\rho_j}[k]\hspace{1mm} \ldots\hspace{1mm}\check{\rho_l}[k]\hspace{1mm} \ldots]^{\top} \in \mbb{R}^{N+N_I+N_O},$$ for which $i\in\Omega$, $j\in\Omega_I$ and $l\in\Omega_O$. Here, $k$ is same as the time-index of the simulation where each time step is of duration $T$. The estimation interval is thus equal to $T$. The equations describing state vector evolution can be be divided into categories:
%\vspace{-0.4cm}

	\vspace{0.1cm}
\begin{subequations}
\label{eq:non_linear_state_evolution}
\begin{itemize}[leftmargin=*]
\item $i\in\Omega\setminus\Omega_I\cup\Omega_O$
\begin{align}
\begin{split}
x_i[k+1] &=  x_i[k]+ \frac{T}{l} \Big(\min\big(\delta_{i-1}[k],\,\sigma_i[k]\big)\\
&\quad-\min\big(\delta_i[k],\,\sigma_{i+1}[k]\big)\Big)\end{split}
\end{align}
\item $i\in\Omega_I$, $j=N+\bar{j}, \bar{j}\in\hat{\Omega}$
\begin{align}
\begin{split}
&x_i[k+1] =  x_i[k]+ \frac{T}{l} \Big(\min\big(\delta_{i-1}[k],\sigma_i[k]\big)\\
&\quad\;-\min\big(\delta_i[k],\,\sigma_{i+1}[k]\big)\\
&\quad\; +\min\big(v_fx_j[k],\xi_{i}[k](\rho_m-x_i[k]),\frac{\xi_i[k]}{w_c}v_f\rho_c\big)\Big)\end{split}
\end{align}
\item  $i\in\Omega_O$, $j=N+N_I+\bar{j}, \bar{j}\in\check{\Omega}$
\begin{align}
\begin{split}
x_i[k+1] &=  x_i[k]+ \frac{T}{l} \Big(\min\big(\delta_{i-1}[k],\,\sigma_i[k])\big)\\
&\quad-\frac{1}{\bar{\beta}_{i}[k]}\min\big(\delta_i[k],\,\sigma_{i+1}[k]\big)\Big)\end{split}
\end{align}
\item $\bar{i}\in\hat{\Omega},i=N+\bar{i},j\in\Omega_I$
\begin{align}
\begin{split}
&x_{i}[k+1] = x_{i}[k]+ \frac{T}{l}\big(\hat{r}_j[k]\\
&\quad-\min\big(v_fx_i[k],\xi_j[k](\rho_m-x_j[k]),\frac{\xi_j[k]}{w_c}v_f\rho_c\big)\big)\end{split}
\end{align}

\item $\bar{i}\in\check{\Omega},i=N+N_I+\bar{i},j\in\Omega_O$
\begin{align}
x_i[k+1] &= x_{i}[k] +\frac{T}{l}\Big(\frac{\beta_{j}[k]}{\bar{\beta}_{j}[k]}\min\big(\delta_j[k],\sigma_{j+1}[k]\big)-\check{s}_{j}[k]\Big)
\end{align}

\end{itemize}
\end{subequations}
%\comment{Suyash: No good connection between the above and the following theorem. Also, we should not having theorems in this paper. Everything should be a proposition. Fix this and clean up the writing, make the next result more compact (it's too long now). You can simply state that the aggregate model (1)--(11) + the above bulleted equations can be written in a compact state space format and give the state-space equations...where A, G, f, Bu are matrices and vector valued functions that are specified in the appendix. Also, where is the proof for this result? You should say that the proof is in Appendix B. Look at previous papers from Sebastian and see how we do this. } 
\vspace{-0.2cm}
\begin{myprs}
%The above equations \eqref{eq:non_linear_state_evolution} describing the evolution
The evolution of traffic density %of a ramp-connected stretched highway system
 described in \eqref{eq:non_linear_state_evolution} can be written in a compact state-space form as follows 
\begin{align}
    \m x[k+1]=\m A\m x[k]+\m G\m f(\m x,\m u)+\m {B_{\mathrm{u}}}\m u[k],
\label{eq:state_space_gen}
\end{align}
\noindent where $\m A \in \mbb{R}^{n\times n}$ for $n := N+N_I+N_O$ represents the linear dynamics of the system, $\m {B_{\mathrm{u}}}\in \mbb{R}^{n\times m}$ for $m := 2+N_I+N_O$ represents the way external inputs affecting the system, $\m f:  \mbb{R}^{n}\times  \mbb{R}^{m}\rightarrow  \mbb{R}^{g}$ is a vector valued function representing nonlinearities in \eqref{eq:non_linear_state_evolution}, and $\m G \in \mbb{R}^{n\times g}$ is a matrix representing the distribution of nonlinearities.\\
%Proof. {\rm See Appendix \ref{apdx:state_space_parameters}}.
\end{myprs}
\vspace{-0.4cm}
The proof and the structures of the matrices and function in \eqref{eq:state_space_gen} are all provided in Appendix \ref{apdx:state_space_parameters}. % of the supplemental document attached with this manuscript submission.
The above result is important as it allows us to perform optimal sensor placement for the traffic dynamical system defined by \eqref{eq:triangular_model}--\eqref{eq:CTM_flow_conservation}, as well as to develop a robust state observer. This also allows for other control-theoretic studies to use this state-space model for other traffic engineering applications. The next section presents our strategy to address the optimal sensor placement problem for the nonlinear traffic dynamical system derived above.

%\comment{Suyash: Does this need to be a proposition? Also everyone knows what a Jacobian structure looks like. Also, why is f uppercase F in the J matrix? Please fix this. Furthermore, I don't think we should have this appendix or this result in the paper. It occupies three pages of appendices and the paper has become excruciatingly long. }
%
%\comment{Suyash: Please re-read the above two sections. You need to  make them more streamlined and easier to read. Also, your appendices need major revision and updates. They read like a noteboook notes, and i understand that that was your first iteration. I'm sure you can make them better :). Also, if there's something tricky in finding the Jacobian since you have the min function and absolute values, then maybe you can mention it here or whatever. Your call. The bottom line: this needs to be presentable.}\\

\section{Observability-based Sensor Placement}\label{sec:traffic_ssp}
In this section we discuss our approach for addressing the traffic sensor placement problem.
The traffic dynamics \eqref{eq:state_space_gen} with measurements can be expressed as
\begin{subequations}\label{eq:state_space_general}
	\begin{align}
	\m x[k+1] &= \mA \m x[k]+\m G\m f(\m x, \m u,k)+ \m {B_{\mathrm{u}}} \m u[k]\label{eq:state_space_general_a}\\
	\m y[k] &= \m \Gamma\mC \m x[k], \label{eq:state_space_general_b}
	\end{align}
\end{subequations}
with $\m f(\m x, \m u,k) := \m f(\m x, \m u)$ where $k$ is used to indicate the time dependence.
In the above model, we introduce $\m y\in\mathbb{R}^p$ to represent the vector of measurements, which corresponds to all of the highway segments equipped with sensors measuring the density.
The matrix $\m C\in\mathbb{R}^{p\times n}$ in \eqref{eq:state_space_general_b} is useful to determine the placement of traffic sensors, whereas $\m \Gamma := \mathrm{Diag}(\m \gamma)$ with $\m \gamma\in\{0,1\}^p$ represents the selection of sensors---that is,
$\gamma_i = 1$ if the $i$-th highway segment is measured (this consequently leads to nonzero row $i$ in $\m C$) and $\gamma_i = 0$ otherwise.
It is assumed throughout the paper that every trajectory of \eqref{eq:state_space_general} lie in the domain of interest $\m \Omega := \mathbfcal{X} \times \mathbfcal{U}$. 
In the context of traffic dynamics, $\mathbfcal{X}$ represents the operating region of the states $\m x$. Since $\m x$ represents the densities, then $\mathbfcal{X} := [0,\rho_m]^n$. Likewise, the set $\mathbfcal{U}$ can be constructed as $\mathbfcal{U} := [0,v_f\rho_c]^m$, which represents the set of admissible inputs.

Having described system \eqref{eq:state_space_general} with sensor placement, we now state the paper's major computational objective: from a given possible set of traffic sensors $\mathcal{G}_{\gamma}$ such that $\m \gamma\in\mathcal{G}_{\gamma}$, find the best (or \textit{optimal}) sensor configuration $\m \gamma^*$ such that system \eqref{eq:state_space_general} is observable, i.e., the system's initial state $\m x_0 := \m x[0]$ can be uniquely determined from a finite set of measurements. 
%\comment{Why is this paragraph and part of the above one here not in the intro? I suggest you wait for Suyash to finish his intro and then try to incorporate the next paragraph in Section I. Do you agree or do you think this should be here?}
In this paper, we opt to formulate the traffic sensor placement using the concept of observability through MHE framework developed in \cite{Haber2017}. {The framework utilizes a series of past measurement data to estimate the initial state at each time window.} The reasons of pursuing this approach are two-fold. First, as argued in \cite{Haber2017}, this approach is considerably more scalable than using empirical observability Gramian and second, we experience numerical issues in applying the method from \cite{Nugroho2019CDC}.
With that in mind, we consider the first $N$ observation window (or $N$ discrete measurements) of system \eqref{eq:state_space_general} expressed in the equation below
\begin{align}
\m h(\m \gamma,\m x_0) := \tilde{\m y} - \m g\left(\m \gamma,\m x_0\right).\label{eq:obs_wind_error}
\end{align}
%\comment{You use $\m w$ twice in the paper to define two different quantities. }
In \eqref{eq:obs_wind_error}, the mapping $\m g:\mbb{N}^p\times\mbb{R}^n\rightarrow\mbb{R}^{Np}$ is given as
\begin{align}
 \hspace{-0.05cm}\m g^{\top}\hspace{-0.07cm}\left(\m \gamma,\m x_0\right) \hspace{-0.07cm}:= \hspace{-0.07cm}\bmat{\left(\m \Gamma\m C \m x[0]\right)^{\top}\,\left(\m \Gamma\m C \m x[1]\right)^{\top}\,\hspace{-0.05cm}\hdots\,\left(\m \Gamma\m C \m x[N-1]\right)^{\top}\hspace{-0.05cm}}\hspace{-0.08cm}.\label{eq:w_function}
\end{align}
Note that $\m g(\cdot)$ defined above is indeed a function of $\m x_0$ since for any $k$ such that $0 < k \leq N-1$ then we have
\begin{align}
\hspace{-0.24cm}\m x[k] = \m A^k\m x[0]+\sum_{j=0}^{k-1}\m A^{k-1-j}\left(\m G\m f(\m x,\m u,j)+\m {B_{\mathrm{u}}} \m u[j]\right). \label{eq:function_x_k}
\end{align}
%\comment{You need to have the time index $j$ inside the summation for $f(\cdot)$}
\vspace{-0.3cm}
\begin{myrem}\label{rem:function_w}
The function $\m g(\cdot)$ in \eqref{eq:w_function} is in fact also dependent on input vector $\m u$, since the $k$-th state $\m x[k]$, as it is seen from \eqref{eq:function_x_k}, relies on $\mathbfcal{U}_s(k-1)$ defined as $$\mathbfcal{U}_s(k-1) := \{\m u[j]\}_{j=0}^{k-1},\;\text{where}\;\m u[j]\in\mathbfcal{U},\;\forall j=0,\hdots,k-1,$$
i.e., the $k$-th sequence of control inputs. The reason why $\m g(\cdot)$ only takes $\m \gamma$ and $\m x_0$ as arguments is due to the assumption that $\mathbfcal{U}_s(i)$ is known for each time index $i$, thus the dependence of $\m g(\cdot)$ towards input vector $\m u$ can be ignored for simplicity.
\end{myrem}
\vspace{-0.1cm}
The term $\tilde{\m y}\in\mbb{R}^{Np}$ in \eqref{eq:obs_wind_error} denotes the stacked $N$ measurements data constructed as
\begin{align*}
\tilde{\m y}^\top := \bmat{\tilde{\m y}[0]^{\top}\;\,\tilde{\m y}[1]^{\top}\;\,\hdots\;\,\tilde{\m y}[N-1]^{\top}}.
\end{align*} 
Obviously, for every initial condition of the system $\m x_0\in\mathbfcal{X}$, it holds that $\m h(\m \gamma,\m x_0) = 0$ such that $\tilde{\m y} = \m g\left(\m \gamma,\m x_0\right)$. Hence, for a fixed $\m \gamma$, the mapping $\m g(\cdot)$ maps the initial state into the $N$ measurements output. The observability of system \eqref{eq:state_space_general} with respect to $\m g(\cdot)$ is formally defined as follows \cite{Hanba2009}.
\vspace{-0.1cm}
\begin{mydef}\label{def:uniform_observability}
The system \eqref{eq:state_space_general} with a prescribed $\m \gamma$ is uniformly observable in $\mathbfcal{X}$ if, for all admissible inputs $\m u[k]\in \mathbfcal{U}$, there exists a finite $N>0$ such that the mapping $\m g\left(\m \gamma,\m x_0\right)$ defined in \eqref{eq:w_function} is injective (one-to-one) with respect to $\m x_0\in \mathbfcal{X}$.
\end{mydef}
\vspace{-0.1cm}
Note that, from Definition \ref{def:uniform_observability}, if $\m g\left(\m \gamma,\m x_0\right)$ is injective with respect to $\m x_0$, then $\m x_0$ can be uniquely determined from the set of measurements $\tilde{\m y}$.
A sufficient condition for $\m g(\cdot)$ to be injective is that the Jacobian of $\m g(\cdot)$ around $\m x_0$, denoted by $\m {J}_w(\cdot)$, is full rank \cite{Hanba2009}. In fact, it is not difficult to show that, if system \eqref{eq:state_space_general} has no nonlinear counterparts, then the Jacobian matrix $\m {J}_w(\cdot)$ reduces to the $N$-step observability matrix for the linear dynamics. 

In this work, we use the concept of observability Gramian to quantify the system's osbervability for a given set of sensors. The Gramian $\m G_v$ (or Gram matrix) for some finite-dimensional real $k$ vectors $\m v_i\in\mbb{R}^{n}$ where $i = 1,2,\hdots,k$ is constructed as $\m G_v := \m V^{\top} \m V$ where $\m V := \bmat{\m v_1 \,|\, \m v_2 \,\hdots\,|\, \m v_k}$. The observability of a discrete-time linear time-invariant system of the form 
	\begin{align}
		\m x[k+1] = \m A \m x[k],\quad  \m y[k] = \m C \m x[k], \label{eq:LTI-DT}
	\end{align}
where $\m x$ is the state and $\m y$ is the output can be determined from the associated Gramian. The corresponding $N$-step observability Gramian is given as $$\mathbfcal{W}_o := \mathbfcal{O}_N^\top \mathbfcal{O}_N = \sum_{i = 0}^{N-1} (\m A^i)^\top \m C^\top \m C \m A^i,$$ where $\mathbfcal{O}_N$ denotes the $N$-step observability matrix. The system \eqref{eq:LTI-DT} is then observable if and only if $\mathbfcal{W}_o$ is of full rank. 
Using the same analogy, the observability Gramian for the nonlinear system \eqref{eq:state_space_general} with respect to $\m \gamma$ around $\m x_0$ can be constructed as
\begin{align}
\m W_o(\m \gamma,\m x_0) := \m {J}_w^{\top}(\m \gamma,\m x_0)\m {J}_w(\m \gamma,\m x_0),\label{eq:obs_gramian_x0}
\end{align} 
where $\m W_o(\cdot)\in\mbb{R}^{n\times n}$ and $\m {J}_w(\cdot)\in\mbb{R}^{Np\times n}$ is given as \cite{Haber2017}
\begin{align}
\m {J}_w(\m \gamma,\m x_0) := \bmat{\m I \otimes \m \Gamma \m C}\times\bmat{\m I\\  \dfrac{\partial \m x[1]}{\partial \m x[0]}\\\vdots\\\dfrac{\partial \m x[N-1]}{\partial \m x[0]}},\label{eq:obs_jacobian_x0}
\end{align}
where for system \eqref{eq:state_space_general}, by applying chain rule, for $0 < k \leq N-1$ the $k$-th partial derivative can be obtained from 
\begin{align*}
\dfrac{\partial \m x[k]}{\partial \m x[0]} = \prod_{j=0}^{k-1} \m A + \m G \dfrac{\partial \m f}{\partial \m x} (\m x,\m u,j),
\end{align*}
where the term $\tfrac{\partial \m f}{\partial \m x}(\cdot)\in\mbb{R}^{g\times n}$ denotes the Jacobian of $\m f(\cdot)$ with respect to $\m x$ evaluated at discrete time $j$. The explicit form of this Jacobian matrix is detailed in Appendix \ref{apdx:calc_jacobian}.
%in the supplemental document attached with this manuscript submission. 
This requires $\m x[j]$ and $\m u[j]$ to be known, which can be obtained from simulating \eqref{eq:state_space_general} for up to $j$-th discrete time. 
In the context of sensor placement problem, the objective is to find the best $\m \gamma$ satisfying $\m \gamma\in\mathcal{G}_{\gamma}$ which maximizes the observability of system \eqref{eq:state_space_general} based on the observability Gramian \eqref{eq:obs_gramian_x0}.   
%To that end, the objective of the sensor placement problem can be rephrased into finding the best $\m \gamma$ satisfying $\m \gamma\in\mathcal{G}_{\gamma}$ which maximizes the observability of system \eqref{eq:state_space_general} around the initial state $\m x_0$. Since $\m x_0$ is unknown a priori, then the observability is studied around an estimate of initial state $\hat{\m x}_0$, which should be sufficiently close to $\m x_0$.

There exist some measures that quantify observability based on Gramian matrix, some of which include the rank, smallest eigenvalue, condition number, trace, and  determinant---see reference \cite{qi2015optimal}. Each measures has a unique characterization vis-à-vis observability. The rank of observability Gramian quantifies the dimension of observable subspace. The smallest eigenvalue of the Gramian is the worst case metric of observability and reflects the difficulty of estimating the actual initial state $\m x_0$ from the given set of measurements $\tilde{\m y}$. That is, the estimation error can be very much sensitive to observation noise when the smallest eigenvalue of the Gramian matrix is relatively close to zero \cite{Krener2009}. The condition number---which measures the ratio between the largest and smallest eigenvalues of the Gramian matrix---on the other hand, accentuates the sensitivity of observability. 

We realize that the Gramian \eqref{eq:obs_gramian_x0} quantifies the observability on a particular initial state $\m x_0$. Large condition number (which implies \textit{ill-conditioned} Gramian matrix) indicates that small change in initial state may produce a significant impact on observability \cite{Krener2009}. The trace of observability Gramian measures the average observability of the system in all directions in state space. Therefore, a larger value of the trace indicates an increase in overall observability \cite{singh2006sensor}. Similar to the trace, determinant of observability Gramian matrix also measures the system's observability in every directions in state space. Nonetheless, as pointed out in \cite{SERPAS2013105}, determinant has the ability to capture information from all elements of the Gramian matrix as well as taking into account information redundancy for the case when multiple sensors are considered. It is also suggested in \cite{qi2015optimal} that the determinant is a better measure for quantifying observability than the trace, since trace tends to overlook (near) zero eigenvalues.             

%\comment{I think here you fail to mention something like: ``Although these metrics have been studied in the literature, two main things remain unclear and rather uninvestigated: \textit{(i)} the computational time performance of the two metrics; \textit{(i)} their impact on the state estimation performance under noise and uncertainty. (maybe you can think of more)." This contribution is important and should also be stated part of the paper's contributions: from now on we should focus on clear takeaways in our papers rather than dumping optimization problems here. The takeaway here from the case studies is clear and should be articulated.}

%Based on the above considerations, 
Based on the above considerations, we settle upon maximizing the determinant and trace of the observability Gramian since, despite their differences discussed previously, both metrics are popular in sensor placement literature through empirical observability Gramian; see  \cite{qi2015optimal,singh2006sensor,SERPAS2013105}, as well as their easy implementation on standard optimization interfaces such as YALMIP \cite{Lofberg2004} and CVX \cite{cvx2020}. The resulting sensor placement problem is given as follows
%can now be stated as follows
\begin{subequations}\label{eq:ssp_gramian_trace}
\begin{align}
(\textbf{P1})\;\;\kappa = \minimize_{\m\gamma} \;\;\;&\begin{cases}
-\mathrm{det}\left(\m W_o(\m \gamma,\hat{\m x}_0)\right), \\
-\mathrm{trace}\left(\m W_o(\m \gamma,\hat{\m x}_0)\right), 
\end{cases} \label{eq:ssp_gramian_trace_1}\\
\subjectto \;\;\;\,& \;\m \gamma\in\mathcal{G}_{\gamma}, \; \m \gamma\in\{0,1\}^p.  \label{eq:ssp_gramian_trace_2}
\end{align}
\end{subequations}
Notice that \textbf{P1} only takes the integer variables $\m \gamma$ as the optimization variable while $\hat{\m x}_0$ is fixed. Since the constraint and objective function are convex (or can be easily transformed into a convex one)\footnote{The trace is linear while determinant is log-concave (and thus logarithmic determinant is concave) on the set of all positive definite matrices \cite{boyd2004cvxbook}.}, then \textbf{P1} is categorized as a convex IP, which can be solved optimally via a BnB algorithm. 
After the observability of the system has been determined from the solutions of \textbf{P1}, then the measurement equation \eqref{eq:state_space_general_b} can be reformulate into $\m y[k] = \tilde{\mC} \m x[k]$ where $\tilde{\mC}\in\mbb{R}^{n\times\tilde{p}}$ is the reduced state-to-output matrix that corresponds to the nonzero rows of $\m\Gamma^*\m C$, where $\m \gamma^*$ is the optimal solution of \textbf{P1}. The traffic density estimation is then performed based on $\m y[k]\in\mbb{R}^{\tilde{p}}$.

\vspace{-0.1cm}
\section{A Robust Observer for State Estimation}\label{sec:observer_design}
The previous section focuses on determining the traffic sensor placement based on the observability of the traffic dynamics through a finite observation window. As such, the obvious approach for performing traffic density estimation is via MHE. However, compared to other approaches such as observers or Kalman filter (KF)-based estimators, MHE for nonlinear systems is generally computationally more expensive since a nonlinear least-square problem has to be solved in every discrete time instance, which make it only suitable for slow dynamical systems or situations where high-powered computational resources are available. To that end, herein we employ a robust observer for estimating traffic density, since observer uses a single gain matrix for all discrete-time instances and only need to be computed once for a given traffic sensors configuration.    

To develop the robust state observer, the nonlinearities in \eqref{eq:state_space_gen} have to be identified. In fact, we prove that these nonlinearities satisfy the Lipschitz continuity condition, which for a vector-valued function is defined as follows.
\vspace{-0.1cm}
\begin{mydef}[Lipschitz Continuity]\label{def:Lipschitz_cont}
	Let $\m f : \mathbb{R}^{n}\times \mathbb{R}^{m} \rightarrow \mathbb{R}^{n}$. Then, $\m f$ is said to be Lipschitz continuous in $\m{\Omega}\subseteq \mathbb{R}^{n}\times \mathbb{R}^{m}$ if there exists a constant $\gamma \in \mathbb{R}_{+}$ such that
	\begin{align}
	\norm{\m f(\m x,\m u)-\m f(\hat{\m x},\m u)}_2 \leq \gamma \norm{\m x -  \hat{\m x} }_2, \label{eq:lipschitz}
	\end{align}
	for all $\m x, \hat{\m x} \in \mathbfcal{X}$ and $\m u \in \mathbfcal{U}$.
\end{mydef} 
\vspace{-0.1cm}
%With that in mind, we present an analytical methodology to derive the Lipschitz constant for the non-linear function $\m f(\cdot)$.
The next proposition shows that the nonlinear mapping $\m f(\cdot)$ in \eqref{eq:state_space_gen} is Lipschitz continuous with Lipschitz constant $\gamma_l$.
\vspace{-0.1cm}
\begin{myprs}\label{prs:lipschitz_constant}
	The non-linear function $\m f(\cdot)$ governing the traffic dynamics \eqref{eq:state_space_gen} and specified in Appendix \ref{apdx:state_space_parameters} is Lipschitz continuous in $\m x$ and $\m u$ with the Lipschitz constant $\gamma_l$ given as
	\begin{align}
	\label{eq:Lipschitz_function_form}
	\gamma_l = h\hspace{0.5mm}(&v_f,\hspace{0.5mm}\omega_c,\hspace{0.5mm}l,\hspace{0.5mm}N, \hspace{0.5mm}N_I,\hspace{0.5mm}N_O,\hspace{0.5mm}\xi_1,\hspace{0.5mm}\xi_2,\hspace{0.5mm}\dots,\hspace{0.5mm}\xi_{N_I},\nonumber\\
	&\hspace{0.5mm}\beta_1,\hspace{0.5mm}\beta_2,\hspace{0.5mm}\dots,\hspace{0.5mm}\beta_{N_O})
	\end{align}
	where the function $h\hspace{0.5mm}(\cdot)$ is defined as in \eqref{eq:Lipschitz_constant}.\\
%	Proof. {\rm See Appendix \ref{apdx:calc_lipschitz_constant}}.
\end{myprs}
\vspace{-0.5cm}
The proof of Proposition \ref{prs:lipschitz_constant} is given in Appendix \ref{apdx:calc_lipschitz_constant}.
% of the supplemental document attached with this manuscript submission.
In what follows, we consider the traffic dynamics \eqref{eq:state_space_general} with the new  measurement matrix $\tilde{\mC}$ and the addition of \textit{unknown inputs}---which may include unmodeled dynamics, disturbance, process noise, and measurement noise---described as 
\begin{subequations}\label{eq:state_space_general_noise}
	\begin{align}
	\m x[k+1] &= \mA \m x[k]+\m G\m f(\m x, \m u)+ \m {B_{\mathrm{u}}} \m u[k]+ \m {B_{\mathrm{w}}} \m w[k] \label{eq:state_space_general_noise_a}\\
	\m y[k] &= \tilde{\mC} \m x[k]+ \m {D_{\mathrm{w}}} \m w[k]. \label{eq:state_space_general_noise_b}
	\end{align}
\end{subequations}
In the above model, the disturbance vector $\m w\in\mathbb{R}^q$ lumps all unknown inputs into a single vector with the corresponding matrices $\m {B_{\mathrm{w}}}$ and $\m {D_{\mathrm{w}}}$ are of appropriate dimensions. These matrices describe how the external disturbances are distributed in the system, where $\m {B_{\mathrm{w}}}$ takes into account any disturbance affecting the system's dynamics while $\m {D_{\mathrm{w}}}$ considers disturbance affecting the measurements data. Both can be constructed from observing empirical traffic data. 
%It is assumed throughout the paper that the nonlinear function $\m f(\cdot)$ in \eqref{eq:state_space_general} is Lipschitz continuous on $\m \Omega$ with the corresponding Lipschitz constant $\gamma_l$ % \comment{wrong notation.}
The observer dynamics for \eqref{eq:state_space_general_noise} are derived from the classical Luenberger observer given as 
\begin{subequations} \label{eq:nonlinear_observer_dynamics}
	\begin{align}	
	\begin{split}{\hat{\m x}}[k+1] &= \mA \hat{\m x}[k+1]+\m G\m f(\hat{\m x}, \m u) \\&\quad+ \m {B_{\mathrm{u}}} \m u[k] + \mL(\m y[k]-\hat{\m y}[k]) 
	\end{split} \\
	\hat{\m y}[k] &= \tilde{\mC} \hat{\m x}[k],
	\end{align}
\end{subequations}
where $\hat{\m x}[k]$ is the state estimation vector, $\hat{\m y}[k]$ is the measurement estimation vector, and $\mL(\m y-\hat{\m y})$ is the Luenberger-type correction term with $\mL\in\mathbb{R}^{n\times \tilde{p}}$. 
By defining the estimation error as $\m e[k]  := \m x [k]- \hat{\m x}[k]$, then from \eqref{eq:state_space_general_noise} and \eqref{eq:nonlinear_observer_dynamics}, the error dynamics can be shown to be in the following form
\begin{subequations} \label{eq:est_error_dynamics}
\begin{align}
	\begin{split}{\m e} [k+1] &= \left(\mA-\mL\tilde{\mC}\right)\hspace{-0.05cm}\m e[k] +\m G \Delta \m f[k]\\&\quad+ \left(\m {B_{\mathrm{w}}}-\mL\m {D_{\mathrm{w}}}\right)\hspace{-0.05cm}\m w[k] \end{split}\label{eq:est_error_dynamics_1} \\
\m z[k] &= \mZ \m e[k],\label{eq:est_error_dynamics_2}
\end{align}
\end{subequations}
where $\Delta \m f[k]:= \m f(\m x, \m u)-\m f(\hat{\m x}, \m u)$ and in particular, $\m z \in \mathbb{R}^{z}$ is the performance output of the error dynamics with respect to the user-defined performance matrix $\mZ\in\mathbb{R}^{z\times n}$.

The goal of the observer is to provide asymptotic estimation error for estimation error dynamics given above. Since the presence of unknown inputs will not make the estimation error to be exactly zero, a robustness metric, referred to as and $\linf$ stability, is employed. The purpose of this metric is to provide numerical assurance on the behavior of performance output $\m z$ against nonzero, time-varying unknown inputs $\m w$.  Our prior work \cite{nugroho2018journal} deals with a robust observer design using the concept of $\linf$ stability for traffic density estimation purpose assuming nonlinear continuous-time traffic dynamics model corresponding with the Greenshield's model. Furthermore, our recent work~\cite{Vishnoi2020} develops a robust $\linf$ observer for the discrete-time model corresponding with the ACTM.  In the sequel, we reproduce a simple numerical procedure from~\cite{Vishnoi2020} to find an observer gain $\m L$ that, if successfully solved,  renders the estimation error dynamics \eqref{eq:est_error_dynamics} to be $\mathcal{L}_{\infty}$ stable with performance level $\mu$. Compared with~\cite{Vishnoi2020}, this paper includes the mathematical proof for the next result. 
\vspace{-0.1cm}
\begin{myprs}\label{thm:l_inf_theorem}
	Consider the nonlinear dynamics \eqref{eq:state_space_general_noise} and observer \eqref{eq:nonlinear_observer_dynamics} in which $\m w \in \mathcal{L}_{\infty}$ and $\m f : \mathbb{R}^{n}\times \mathbb{R}^{m} \rightarrow \mathbb{R}^{n}$ is locally Lipschitz in $\m \Omega$ with Lipschitz constant $\gamma_l$. If there exist $\mP\in\mathbb{S}^n_{++}$, $\mY\in\mathbb{R}^{n\times p}$, $\epsilon,\mu_0,\mu_1,\mu_2\in \mathbb{R}_{+}$, and $\alpha \in \mathbb{R}_{++}$ such that the following optimization problem is solved
	%\vspace{-0.5cm}
	%	\begin{mdframed}[style=MyFrame]
	\begin{subequations}\label{eq:l_inf_theorem2}
		\vspace{-0.0cm}
		\begin{align}
		%	&\mathcal{L}_{\infty}-\mathrm{Observer} \nonumber\\
		%\mc{L}_{\infty}\text{-}\mathrm{BMI}:
		&(\mathbf{P2})\;\;\minimize_{\m P, \m Y, \epsilon, \alpha, \mu_{0,1,2}} \quad \mu_0\mu_1 + \mu_2 \label{eq:l_inf_theorem2_0}\\
		&\subjectto \;\;\nonumber \\
		&\bmat{  (\alpha-1)\mP+\epsilon\gamma_l^2\mI&*&*&*\\
			\m O&-\epsilon\mI &*&*\\
			\mO & \m O & -\alpha\mu_0\mI &*\\
			\mP\mA-\mY\mC&\m P \m G& \mP\m {B_{\mathrm{w}}}-\mY\m {D_{\mathrm{w}}} &-\mP} \preceq 0 \label{eq:l_inf_theorem2_1}\\
		&\bmat{-\mP & * & * \\
			\mO & -\mu_2\mI & *\\
			\mZ & \mO & -\mu_1\mI}\preceq 0, \label{eq:l_inf_theorem2_2}
		\end{align}
	\end{subequations}
	%	\end{mdframed}
	\noindent 
	then the error dynamics given in \eqref{eq:est_error_dynamics} is $\mathcal{L}_{\infty}$ stable with performance level $\mu = \sqrt{\mu_0\mu_1+\mu_2}$ for performance output given as $\m z = \mZ \m e$ and observer gain $\mL$ computed as $\mL = \mP^{-1}\mY$.
\end{myprs}   
The proof of Proposition \ref{thm:l_inf_theorem} is provided in Appendix \ref{apdx:obs_proof}.
% of the supplemental document attached with this manuscript submission.
Note that optimization problem described in \eqref{eq:l_inf_theorem2} is nonconvex due to bilinearity appearing in the form of $(\alpha-1)\mP$, $\alpha\mu_0$, and $\mu_0\mu_1$. To get a convex one, one can simply fix the values of $\alpha$ and $\mu_0$ or $\mu_1$, and solve for the other variables using any SDP solvers.  
%Second, although SDPs are known to be not scalable for large-scale problems, state-space matrices $\m A, \m C, \ldots$ characterizing the traffic density evolution are extremely sparse, thereby enabling scalable computation of the observer gain $\m L$ in Theorem~1.
Having described our strategies for placing traffic sensors and estimating traffic density via a robust observer framework, we then perform numerical tests of these strategies through various case studies---presented in the ensuing sections.

%\comment{This section on the observer design can be shorter. We should not have Definition or the Remark (I mean Remark 2 isn't a remark anyway). We also don't cite our ACC paper here for some reason. You also suddenly have $B_w \m w$ in the dynamics without explaining how we obtain that, and so on. This section needs some re-envisioning. }
\section{Case Study: Results and Analysis}\label{sec:cases_tudy}

%\subsection{Numerical Simulation Setup}\label{ssec:simulation_setup}
This section demonstrates the proposed approach for determining traffic sensors location and estimating traffic density through $\linf$ observer. 
%\comment{You need to add an itemized list of computational research questions that we're investigating in this case study section. There's four-five high-level computational questions. List them here. In fact, \textbf{every single paragraph} in this case study section reveals an interesting observation. Either list these questions here or maybe have the research questions + their answers in a table towards the end of the case studies section. I think this would be so nice.} 
Specifically, in this numerical study we attempt to answer the following questions:
\begin{enumerate}
	\item[--] \textit{Q1}: How does the observation window length relate to traffic observability and initial state estimation?
	 \item[--] \textit{Q2}: How computationally efficient are the determinant and trace observability metrics in terms of solving the sensor placement problem?
	 \item[--] \textit{Q3}: How do actual and presumed initial states affect the resulting sensors' location? Are sensor placements robust to changes in initial states?
	 \item[--] \textit{Q4}: How reliable is the state estimation when using sensor locations obtained from utilizing determinant and trace observability metrics?
	 \item[--] \textit{Q5}: How does the theory-driven placement of sensors compare to randomized and uniform sensor placement strategies?
\end{enumerate}
All simulations are performed using MATLAB R2019a running on 64-bit Windows 10 with 3.4GHz Intel\textsuperscript{R} Core\textsuperscript{TM} i7-6700 CPU and 16 GB of RAM with YALMIP \cite{Lofberg2004} as the interface to solve all IP and convex SDP. Throughout the section, all highways are configured with $v_f = 28.8889$ m/s ($65$ mph), $w_c = 6.6667$ m/s ($15$ mph), $\rho_c = 0.0249$ vehicles/m ($40$ vehicles/mile), $\rho_m = 0.1333$ vehicles/m, and $l = 400$ m. The discrete-time step is chosen to be $T = 1$ \textrm{sec}. 
\begin{figure*}
	\vspace{-0.0cm}
	\centering 
	\subfloat[\label{fig:scenario_1a}]{\includegraphics[keepaspectratio=true,scale=0.610]{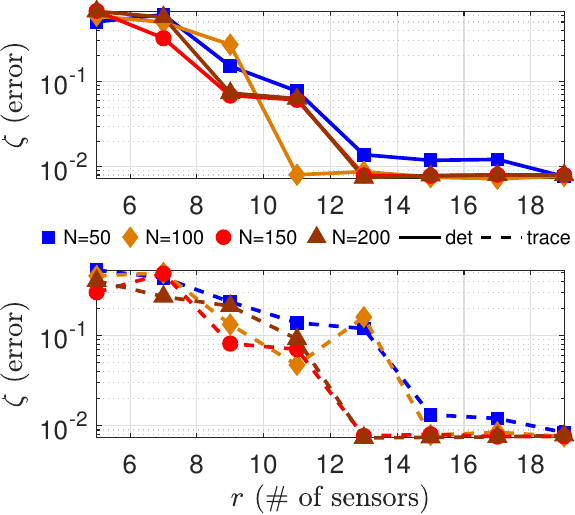}}{}{}\vspace{-0.0cm}
	\subfloat[\label{fig:scenario_1b}]{\includegraphics[keepaspectratio=true,scale=0.610]{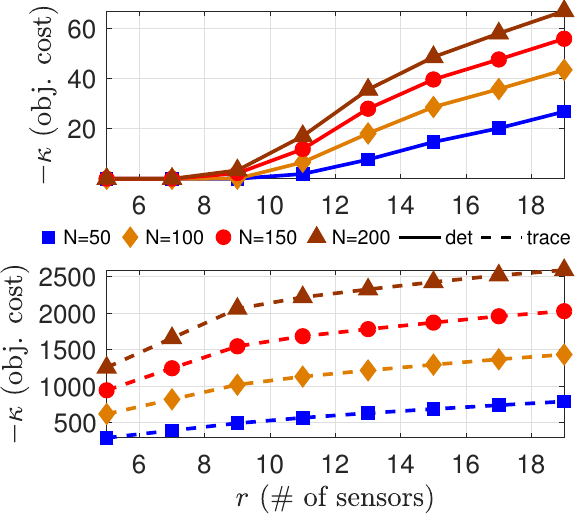}}{}{}\hspace{-0.0cm}
	\subfloat[\label{fig:scenario_1c}]{\includegraphics[keepaspectratio=true,scale=0.610]{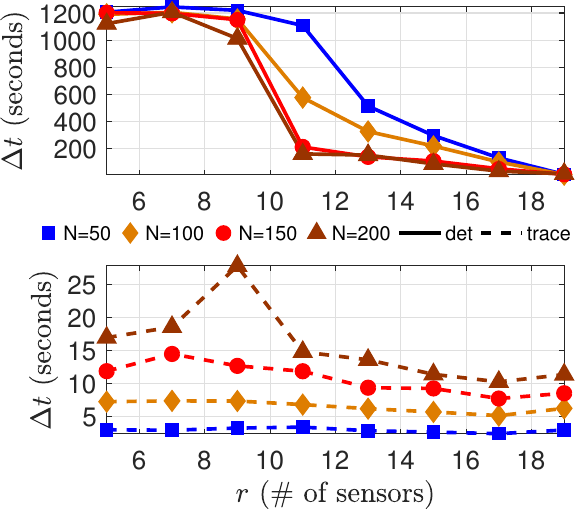}}{}{}\hspace{-0.0cm}\vspace{-0.1cm}
	\caption{Results on observability analysis for different number of sensor allocations: (a) the relative error $\zeta$, (b) inverse optimal value of \textbf{P1}, denoted by $-\kappa$, and (c) total computational time, which includes the overall time spent for solving \textbf{P1} and \textbf{P3}.}
	\label{fig:scenario_1}\vspace{-0.1cm}
\end{figure*} 

\begin{figure}[t]
	%\vspace{-0.2cm}
	\centering 
	\subfloat[\label{fig:scenario_2}]{\includegraphics[keepaspectratio=true,scale=0.61]{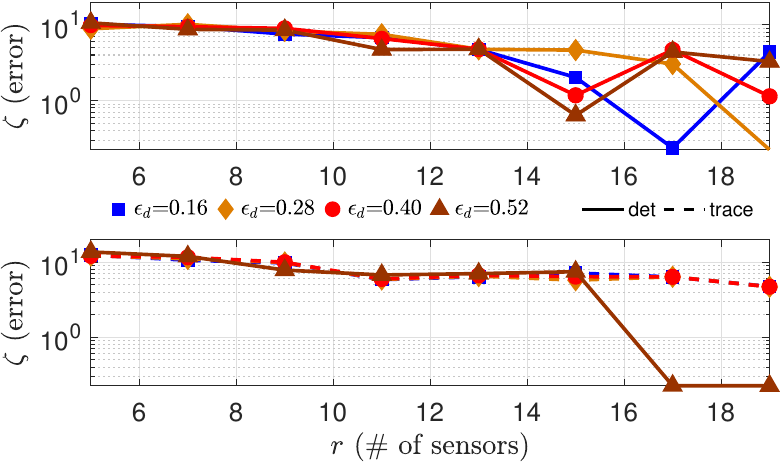}}{}\vspace{-0.1cm}\vspace{-0.2cm}
	\subfloat[\label{fig:scenario_3}]{\includegraphics[keepaspectratio=true,scale=0.61]{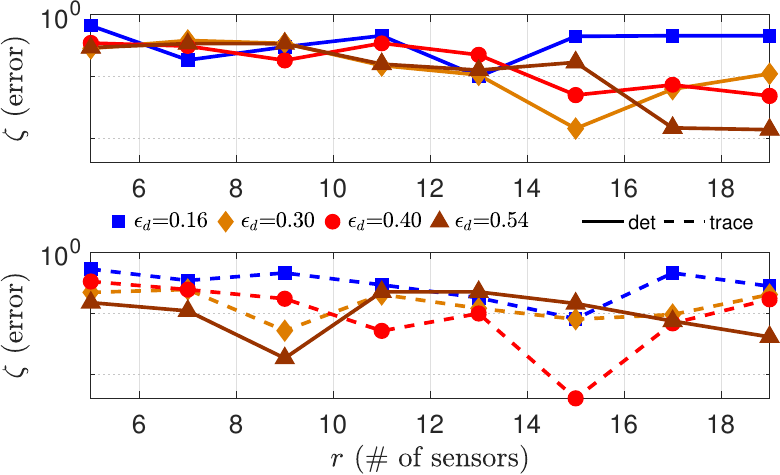}}{}\vspace{-0.1cm}\hspace{-0.0cm}\vspace{-0.02cm}
	\caption{Relative error $\zeta$ for (a) various presumed initial state $\hat{\m x}_0$ and (b) actual initial state ${\m x}_0$.}
	\label{fig:scenario_2_and_3}\vspace{-0.3cm}
\end{figure}

\vspace{-0.3cm}
\subsection{Observability Analysis for Traffic Sensor Placement}\label{ssec:obs_sensor_allocation}
Herein, we perform a numerical analysis on sensor placement approach through traffic network's observability discussed in Section \ref{sec:traffic_ssp}. The highway, referred to the rest of the section as \textit{Highway A}, stretches for approximately $3.2$ miles with $13$ segments on the mainline, $4$ on-ramps, and $4$ off-ramps such that there are $21$ segments in total. The objective of traffic sensor placement problem translates to finding the set of $r$ highway segments that must be equipped with traffic sensors such that the entire highway traffic is observable, which is carried out by solving \textbf{P1}. Realize that \textbf{P1} is classified as a convex \textit{integer programming} (IP) since $\hat{\m x}_0$, the presumed initial state, is fixed. \textbf{P1} is solved using YALMIP's branch-and-bound (BnB) algorithm \cite{Lofberg2004} along with MOSEK \cite{Andersen2000} solver. 

%\comment{Problem \textbf{P1} is redundant. \textbf{P1} alone is sufficient. Fix this. }

In the first instance of observability analysis, problem \textbf{P1} is solved with fixed actual initial state ${\m x}_0$ and presumed initial state $\hat{\m x}_0$, which are generated randomly in $(0,\rho_m]^n$. \textbf{P1} is then solved with different observation windows $N$ for varying number of allocated sensors $r$. In this case, $r = \lceil p\times r_p\%\rceil$ where $\lceil\cdot\rceil$ is the ceiling function, $p = n = 21$, and $r_p\%$ represents the percentage of sensor's allocation (from $20\%$ to $90\%$). The number of allocated sensors is a convex constraint such that, with respect to \eqref{eq:ssp_gramian_trace_2}, $\m \gamma\in\mathcal{G}_{\gamma}$ is equivalent to $\sum_i \gamma_i = r$. 

We solve \textbf{P1} for both observability metrics---the trace and determined metrics discussed in Section~\ref{sec:traffic_ssp}. It is suggested in \cite{qi2015optimal} to use logarithmic determinant as opposed to pure determinant to avoid numerical problems. However, since YALMIP does not have option for realizing logarithmic determinant objective function particularly for IP, then \textbf{P1} is solved using \texttt{geomean} function, which is equivalent to $\sqrt[n]{\mathrm{det}(\cdot)}$. The maximum iteration of YALMIP's BnB algorithm is chosen to be $10^4$. If the optimality gap is not satisfied with the default tolerance after the maximum iteration is achieved, the resulting solution is taken as the final outcome.

In this numerical study, we put our interest on comparing the relative error $\zeta$, the inverse optimal value of \textbf{P1}, denoted by $-\kappa$, and total computational time $\Delta t$. The relative error $\zeta$ represents the difference between actual initial state ${\m x}_0$ and the \textit{computed} initial state $\tilde{\m x}_0$, which is obtained from solving the following nonlinear least-square problem \cite{Haber2017}
\begin{subequations}\label{eq:nonlinear_least_square}
	\begin{align}
	(\textbf{P3})\;\;\minimize_{\tilde{\m x}_0} \;\; &\norm{\tilde{\m y}_{\gamma} - \m g_{\gamma}\left(\tilde{\m x}_0\right)}^2_2 \label{eq:nonlinear_least_square_1}\\
	\subjectto \;\;\;& \;\m 0 \leq \tilde{\m x}_0\leq \m 1\times\rho_m.  \label{eq:nonlinear_least_square_2}
	\end{align}
\end{subequations}
In \eqref{eq:nonlinear_least_square_1}, $\tilde{\m y}_{\gamma}$ is obtained from simulating the undisturbed traffic dynamics \eqref{eq:state_space_general} with a prescribed set of sensors $\m\gamma$ from $k = 0$ to $k = N$ and initial state ${\m x}_0$. Likewise, the nonlinear mapping $\m g_{\gamma}(\cdot)$ is constructed as in \eqref{eq:w_function} with the same sensor combination $\m\gamma$. Problem \textbf{P3} is solved by using the MATLAB data-fitting function \texttt{lsqnonlin}, which implements a \textit{trust region reflective} algorithm \cite{Coleman1996}. The optimality tolerance is set to be $10^{-6}$ and an initial guess equals to $\hat{\m x}_0$. After $\tilde{\m x}_0$ is obtained as the solution of \textbf{P3}, $\zeta$ is computed as
\begin{align*}
\zeta :=\dfrac{\norm{\tilde{\m x}_0-{\m x}_0}_2}{\norm {{\m x}_0}_2}.
\end{align*}   
Note that $\zeta$ appraises the quality of sensor placement, as smaller $\zeta$ suggests that the given sensor configuration is capable to provide a more accurate estimate of the initial state. Of course, this is reasonable if \textbf{P3} returns an \textit{optimal} solution. Nonetheless, as \textbf{P3} is a nonconvex optimization problem, then it is highly unlikely that optimal solutions can be obtained. To reduce the chance of variability, the actual ${\m x}_0$ and presumed $\hat{\m x}_0$ initial states are fixed for this particular test. %The values for ${\m x}_0$ and $\hat{\m x}_0$ are random

\setlength{\textfloatsep}{10pt}
\begin{table*}
	%	\footnotesize
	\scriptsize
	\vspace{-0.31cm}%\hspace{-0.21cm}
	\centering 
	\caption{Sensor locations for two different observation lengths and observability metrics. The notation $(x)+\{y\}$ means that the set of sensors $\{y\}$ is appended to that of row $(x)$ from the same column. The italicized numbers indicate suboptimal sensor locations.}
	\label{tab:sensor_loc_scenario_1}
	\vspace{-0.2cm}
	\renewcommand{\arraystretch}{1.5}
	\begin{tabular}{l|l|l|l|l}
		\midrule \hline
		\multirow{2}{*}{$r$} & \multicolumn{2}{c|}{$\mathrm{det}$}                             & \multicolumn{2}{c}{$\mathrm{trace}$}                           \\ \cline{2-5} 
		& \multicolumn{1}{c|}{$N = 100$} & \multicolumn{1}{c|}{$N = 200$} & \multicolumn{1}{c|}{$N = 100$} & \multicolumn{1}{c}{$N = 200$} \\ \hline
		$5$	&       $\{1,6,10,16,20\}$           &   $\{1,6,10,14,15\}$               &        $\{6,14,17,19,20\}$           &   $\{6,14,17,19,20\}$                \\ \hline
		$7$	&     $\mathit{\{1,6,15,16,17,20,21\}}$             &     $\mathit{\{1,6,14,16,17,20,21\}}$             &    $(5)+\{16,21\}$               &    $(5)+\{15,16\}$               \\ \hline
		$9$	&  $(7)+\mathit{\{14,19\}}$                &   $\{1,14,15,16,17,18,19,20,21\}$               &    $(7)+\{15,18\}$               &   $(7)+\{18,21\}$                \\ \hline
		$11$	&     $(9)+\{10,18\}$             &      $(9)+\{6,10\}$            &      $(9)+\{8,10\}$          &    $(9)+\{8,10\}$              \\ \hline
		$13$	&    $(11)+\{8,12\}$              &    $(11)+\{8,12\}$              &     $(11)+\{11,13\}$              &      $(11)+\{1,11\}$             \\ \hline
		$15$	&   $(13)+\{2,5\}$               &    $(13)+\{2,5\}$              &        $(13)+\{1,2\}$           &     $(13)+\{7,9\}$              \\ \hline
		$17$	&      $(15)+\{7,11\}$            &       $(15)+\{9,11\}$           &       $(15)+\{7,9\}$            &       $(15)+\{2,13\}$            \\ \hline
		$19$	&      $(17)+\{3,13\}$             &     $(17)+\{3,13\}$             &  $(17)+\{3,12\}$                 &   $(17)+\{3,12\}$                \\ 
		\toprule \bottomrule
	\end{tabular}
	\vspace{-0.1cm}
\end{table*}
\setlength{\floatsep}{10pt}  

The results of this numerical study are depicted in Fig. \ref{fig:scenario_1}. In particular, it can be seen from Fig. \ref{fig:scenario_1a} that larger observation window yields smaller relative error. Notice that this behavior is reflected better from the trace but less evident from the determinant---the variations are presumed to be caused by suboptimal solutions obtained from solving \textbf{P3}. Fig. \ref{fig:scenario_1b} also indicates this behavior: larger observation window yields better optimal value \textbf{P1}, which in turn implies that the observability measure is directly proportional with observation window. This corroborates first principles in control theory; more available data (i.e., a larger observation window) results in a better system observability. The corresponding total computational time for both \textbf{P1} and \textbf{P3} are reported in Fig. \ref{fig:scenario_1c}. It can be observed that, in this test where YALMIP's BnB algorithm is utilized, employing the trace in \textbf{P1} is more efficient than doing so with the determinant. 

\setlength{\textfloatsep}{10pt}
\begin{table*}
	%	\footnotesize
	\scriptsize
	\vspace{-0.0cm}%\hspace{-0.21cm}
	\centering 
	\caption{Sensor's locations for two different Euclidean distances between actual ${\m x}_0$ and presumed $\hat{\m x}_0$ initial states where ${\m x}_0$ is fixed and $\hat{\m x}_0$ varies. The notation $(x)+\{y\}$ means that the set of sensors $\{y\}$ is appended to that of row $(x)$ from the same column. The italicized numbers indicate suboptimal sensor locations.}
	\label{tab:sensor_loc_scenario_2_and_3}
	\vspace{-0.2cm}
	\renewcommand{\arraystretch}{1.5}
	\begin{tabular}{l|l|l|l|l}
		\midrule \hline
		\multirow{2}{*}{$r$} & \multicolumn{2}{c|}{$\mathrm{det}$}                             & \multicolumn{2}{c}{$\mathrm{trace}$}                           \\ \cline{2-5} 
		& \multicolumn{1}{c|}{$\epsilon_d = 0.16$} & \multicolumn{1}{c|}{$\epsilon_d = 0.40$} & \multicolumn{1}{c|}{$\epsilon_d = 0.16$} & \multicolumn{1}{c}{$\epsilon_d = 0.40$} \\ \hline
		$5$	&       $\mathit{\{1,6,15,17,20\}}$           &   $\mathit{\{1,9,15,16,19\}}$               &        $\{14,15,18,19,20\}$           &   $\{14,15,18,19,20\}$                 \\ \hline
		$7$	&     $\mathit{\{1,5,16,17,18,20,21\}}$             &     $\mathit{\{1,8,15,17,18,19,21\}}$             &    $(5)+\{16,17\}$               &    $(5)+\{16,17\}$               \\ \hline
		$9$	&  $\mathit{\{1,6,11,14,15,16,18,19,20\}}$                &   $\mathit{\{1,5,10,14,16,17,18,19,21\}}$               &    $(7)+\{3,21\}$               &   $(7)+\{3,21\}$                \\ \hline
		$11$	&     $(7)+\{9,14,15,19\}$             &      $(5)+\{5,14,17,18,20,21\}$            &      $(9)+\{6,9\}$          &    $(9)+\{6,9\}$              \\ \hline
		$13$	&    $(9)+\{3,9,17,21\}$              &    $(5)+\{3,6,11,14,17,18,20,21\}$              &     $(11)+\{2,5\}$              &      $(11)+\{2,5\}$             \\ \hline
		$15$	&   $(11)+\{3,7,11,12\}$               &    $(9)+\{3,6,8,12,15,20\}$              &        $(13)+\{1,8\}$           &     $(13)+\{1,8\}$              \\ \hline
		$17$	&      $(13)+\{2,5,8,12\}$            &       $(13)+\{2,5,8,12\}$           &       $(15)+\{4,7\}$            &       $(15)+\{4,7\}$            \\ \hline
		$19$	&      $(17)+\{4,10\}$             &     $(17)+\{4,10\}$             &  $(17)+\{10,11\}$               &   $(17)+\{10,11\}$                \\ 
		\toprule \bottomrule
	\end{tabular}
	\vspace{-0.2cm}
\end{table*}
\setlength{\floatsep}{10pt}

\begin{figure}[t]
	%\vspace{-0.2cm}
	\centering 
	{\includegraphics[keepaspectratio=true,scale=0.61]{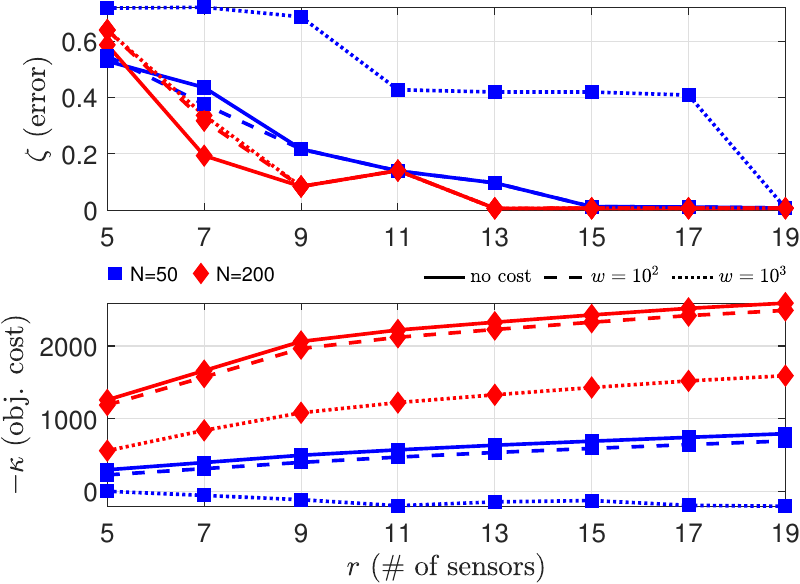}}{}\vspace{-0.1cm}\vspace{-0.2cm}
	\caption{Comparison of relative error (top) and the trace of the Gramian (bottom) between uniform and different (randomized) sensor costs.}
	\label{fig:scenario_9_hwy1}\vspace{-0.3cm}
\end{figure}

The resulting sensor's location for observation window $N = 100$ and $N = 200$ are shown in Tab. \ref{tab:sensor_loc_scenario_1}. For the same number of sensors and observation windows, the sensor locations for determinant and trace objective functions are quite distinct. This observation is expected as both metrics measure different aspect on observability Gramian. Interestingly, as it is seen from Tab. \ref{tab:sensor_loc_scenario_1}, the results from using trace as observability metric suggest a \textit{modular-like behavior}\footnote{A set function $f: 2^{\mathcal{V}}\rightarrow \mbb{R}$, in which $\mathcal{V} = \{1,2,\hdots,p\}$ for $p \in \mbb{N}$, is modular if and only if for any $\mathcal{S}\subseteq\mathcal{V}$ it can be expressed as $f(\mathcal{S}) = w(\emptyset)+\sum_{s\in \mathcal{S}} w(s)$ for some weight function $w:\mathcal{V}\rightarrow \mbb{R}$ \cite{Lovasz1983}.}:
 the sensor combination from a particular row covers the ones from its previous rows. 
%\comment{If someone in transportation is reading this, they won't understand what you mean. Maybe a short intro to modularity and the role it plays is important here (maybe a footnote? be brief} 
This modularity property is also observed from the determinant objective function, despite some irregularities for the first three rows, which is caused by suboptimal solutions. This phenomena is in tandem with the conjecture stated in \cite{Haber2017}, claiming that the logarithmic determinant objective function in \textbf{P1} is submodular.    

Next, we asses the variability of sensor placement due to the differences in actual ${\m x}_0$ and presumed $\hat{\m x}_0$ initial states. The numerical experiment is carried out as follows. First, we set ${\m x}_0$ to be fixed and use different values for $\hat{\m x}_0$ which are randomly generated inside $(0,\rho_m]^n$ in such a way that their Euclidean distances, computed as $\epsilon_d := \norm{\hat{\m x}_0-{\m x}_0}_2$, are unique to each other. For this scenario, the results are given in Fig. \ref{fig:scenario_2}, which shows the resulting relative error. In another scenario, $\hat{\m x}_0$ is fixed while ${\m x}_0$ varies---see Fig. \ref{fig:scenario_3} for the results. 
%Both figures indicate that, despite of some variations, the relative errors $\zeta$ seem to be less influenced by ${\m x}_0$ and $\hat{\m x}_0$. 
It can be seen from both figures that, in general, varying $\hat{\m x}_0$ while fixing ${\m x}_0$ yields higher relative error magnitudes than fixing $\hat{\m x}_0$ and varying ${\m x}_0$.
However, when ${\m x}_0$ is fixed, $\zeta$ experiences minimal variations at least for $r < 14$, whereas more variations occur on the other scenario. 
Again, these variations are most likely attributed to the suboptimal solutions obtained from solving \textbf{P3} via \texttt{lsqnonlin}. It is also observed that larger distance $\epsilon_d$ does not necessarily result in larger relative error $\zeta$, suggesting that the proposed method is rather resilient towards the values of actual and presumed initial states.  

The corresponding sensor combinations for scenario when ${\m x}_0$ is fixed and  $\hat{\m x}_0$ varies with two different Euclidean distances are provided in Tab. \ref{tab:sensor_loc_scenario_2_and_3}. It can be seen that, for the trace, the sensor combinations are the same, indicating that it is robust to different actual and presumed initial states. The same pattern is also observed for the determinant, despite being less apparent. The first three rows are expected since they are suboptimal, hence the different sensor combinations. However, when the solutions are optimal, the resulting sensor combinations are similar except for $r = 15$. Similar results are also obtained for scenario when ${\m x}_0$ varies and  $\hat{\m x}_0$ is fixed. %The resulting sensor combinations are not shown for brevity.

Lastly, in this section we investigate the impact of different sensor costs towards the quality of sensor's configuration. To that end, the following problem is solved
\begin{subequations}\label{eq:ssp_trace_cost}
\begin{align}
\hspace{-0.42cm} (\textbf{P4}) \; \; \, \minimize_{\m\gamma,\, t \geq 0} \;\;\;&
	-\mathrm{trace}\left(\m W_o(\m \gamma,\hat{\m x}_0)\right) + w\times t \\
\hspace{-0.5cm}	\subjectto \;\; & \;\m \gamma\in\mathcal{G}_{\gamma}, \; \m \gamma\in\{0,1\}^p, \;\m c^\top \m \gamma \leq t.
\end{align}
\end{subequations}
Realize that \textbf{P4} is \textbf{P1} with the addition of sensor costs and only considers the trace function since YALMIP faces numerical errors while \texttt{geomean} function is fused with sensor costs. 
 In this problem, $w \in \mbb{R}_{++}$ is used to determine the priority between the two objectives and $\m c\in \mbb{R}_{+}^{p}$ represents the costs for all sensors. Two different values for $w$ are used, where $w \geq 10^2$ to signify the sensor costs. The sensor costs are generated randomly such that $c_i \in [0,1]$. The results are shown in Fig. \ref{fig:scenario_9_hwy1}, from which it can be seen that the relative errors and the degrees of observability for sensor locations determined by considering sensor costs are considerably worse than not using sensor costs at all. This happens because \textbf{P4} has to simultaneously maximizing system's observability and minimizing the sensor costs, thereby giving sensor locations that have lower degree of observability.

%\subsection{Comparing Branch-and-Bound with Greedy Heuristics}\label{ssec:obs_sensor_greedy}
%The results from Section \ref{ssec:obs_sensor_allocation} empirically indicate the submodularity property of \textbf{P1}, which is also conjectured in \cite{Haber2017} for logarithmic determinant maximization. As it can be viewed from Fig. \ref{fig:scenario_1c}, solving \textbf{P1} with determinant objective function indeed requires a tremendous computational effort as opposed to using the trace. To combat this undesirable result, assuming that maximizing the determinant of Gramian matrix \eqref{eq:obs_gramian_x0} is submodular, we can utilize greedy heuristics---in lieu of the classical BnB algorithm---to solve \textbf{P1} more efficiently.   

\begin{figure}[t]
	%\vspace{-0.2cm}
	\centering 
	\subfloat[\label{fig:scenario_5_error}]{\includegraphics[keepaspectratio=true,scale=0.61]{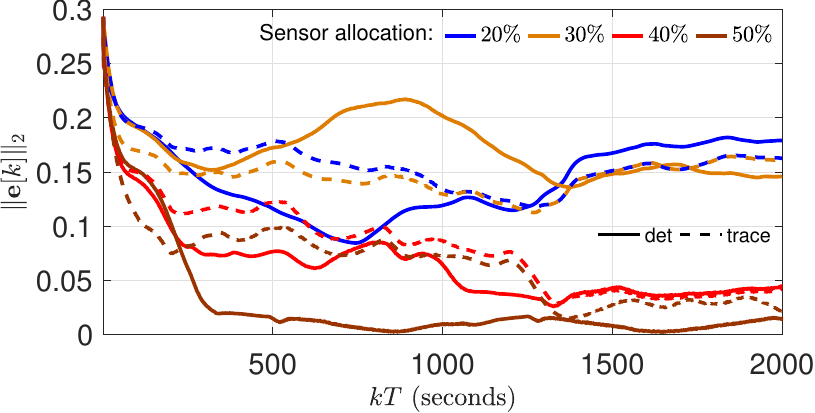}}{}\vspace{-0.1cm}\vspace{-0.2cm}
	\subfloat[\label{fig:scenario_5_rmse}]{\includegraphics[keepaspectratio=true,scale=0.61]{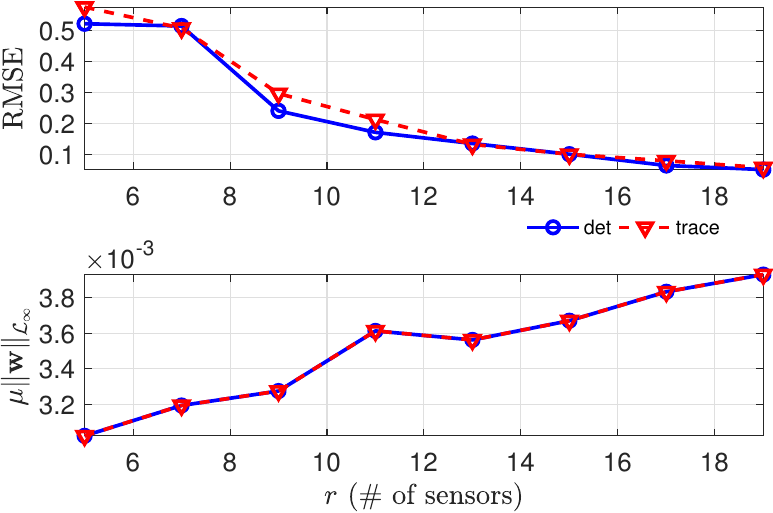}}{}\vspace{-0.1cm}\hspace{-0.0cm}\vspace{-0.02cm}
	\caption{State estimation results via $\linf$ observer for different allocations of sensor: (a) the estimation error and (b) the estimation error performance bounds.}
	\label{fig:scenario_5}\vspace{-0.3cm}
\end{figure}

\begin{figure}[t]
	%\vspace{-0.2cm}
	\centering 
	\subfloat[\label{fig:scenario_8_rmse_hwy1}]{\includegraphics[keepaspectratio=true,scale=0.61]{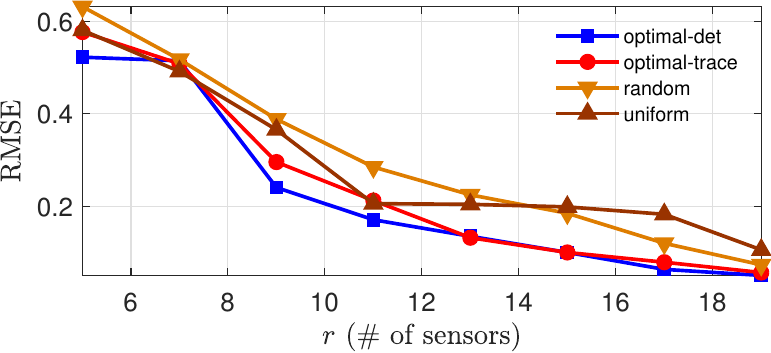}}{}\vspace{-0.1cm}\vspace{-0.2cm}
	\subfloat[\label{fig:scenario_7_rmse_hwy2}]{\includegraphics[keepaspectratio=true,scale=0.61]{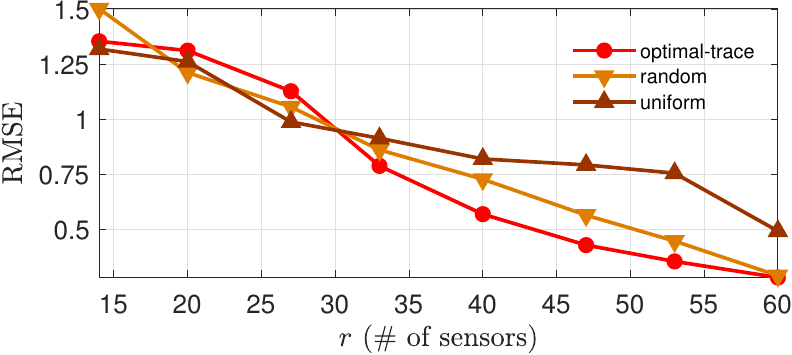}}{}\vspace{-0.1cm}\hspace{-0.0cm}\vspace{-0.02cm}
	\caption{RMSE for optimal and randomized sensor placement for (a) Highway A and (b) Highway B. The simulations with randomized sensor placement are performed $10$ times and the presented results correspond to the average values.}
	\label{fig:scenario_7_8_rmse}\vspace{-0.1cm}
\end{figure}

\vspace{-0.2cm}
\subsection{Traffic Density Estimation with Various Sensor Allocations}\label{ssec:traffic_estimation}
This section investigates different traffic sensor placements---obtained from the previous section---for robust traffic density estimation purpose. To that end, we implement the $\linf$ observer proposed in Section \ref{sec:observer_design} and set $\alpha = 0.1$ and $\mu_1 = 10^4$ to get a convex optimization problem out of \textbf{P2}, in which SDPNAL+ \cite{Yang2015} is used to solve the problem since MOSEK returns numerical problems. The performance matrix is chosen to be $\mZ = 0.01\mI$. 
Herein we simulate Highway A with final time 
$k_f = 2000$ and invoke Gaussian noise with covariance matrices $Q = \nu\mI$ and $R = \nu\mI$ where $\nu = 10^{-3}$ with corresponding unknown input matrices 
$\m {B_{\mathrm{w}}} = \bmat{\m {B_{\mathrm{u}}}\;\;\m O}$ and $\m {D_{\mathrm{w}}} = \bmat{\m O \;\; \m I}$, which simulate process and measurement noise. 
%The simulations are performed on Highway A with final time 
%$k_f = 2000$.
% (or equivalently $t_f =k_fT = 2000$ sec).

In this part of numerical study, 
%we perform traffic density estimation with various number of sensors. In particular, 
we are focused on finding out whether different observability metrics used to solve \textbf{P1} have explicit impact on the quality of traffic density estimation, since according to Tab. \ref{tab:sensor_loc_scenario_1}, different observability metrics return distinctive sensor locations. With that in mind, we compare the resulting estimation errors $\m e[k]$, root-mean-square error (RMSE)---which is computed from
\begin{align*}
\mathrm{RMSE} &= \sum_{i=1}^{n} \sqrt{\frac{1}{k_f}\sum_{k=0}^{k_f}(e_i[k])^2},
\end{align*} 
and estimation error performance bounds $\mu\norm {\m w}_{\linf}$, which is obtained from multiplying the worst case disturbance $\norm {\m w}_{\linf}$ with performance index $\mu$ retrieved from solving \textbf{P2}. The sensor locations used in this test, for both observability metrics, are collected from Tab. \ref{tab:sensor_loc_scenario_1} that correspond to observation window $N = 200$. 

The results of this numerical test are illustrated in Fig. \ref{fig:scenario_5}.
It is indicated from Fig. \ref{fig:scenario_5_rmse} that, as more sensors are utilized, the estimation error is decreasing. However, it is of importance to note that using determinant as the observability metric yields better state estimation results %\comment{better sensor locations? you mean better state estimation corresponding with such placement??} 
than doing so with trace, since the RMSE for determinant is slightly smaller. This showcases the prominence of determinant as opposed to trace to quantify observability. Interestingly, the estimation error performance bounds are increasing as more sensors are utilized with no noticeable difference for determinant and trace. This suggests that, regardless of sensor locations, adding sensors increases the susceptibility of the entire network towards measurement noise. The trajectories of estimation error for different number of sensors are given in Fig. \ref{fig:scenario_5_error}, from which it can be observed that utilizing more sensors yields smaller and faster convergence on estimation error---notice that this observation is in accordance with RMSE shown in Fig. \ref{fig:scenario_5_rmse}. It is also observed from Fig. \ref{fig:scenario_5_error} that, in general, the determinant returns smaller estimation error compared to the trace. A slight exception occurs when the sensor allocation is at $30\%$ (also see Fig. \ref{fig:scenario_5_rmse}). This is suspected to be caused by suboptimal sensor placement---see Tab. \ref{tab:sensor_loc_scenario_1}.

\begin{figure}[t]
	%\vspace{-0.2cm}
	\centering 
	\subfloat[\label{fig:scenario_8_performance_hwy1}]{\includegraphics[keepaspectratio=true,scale=0.61]{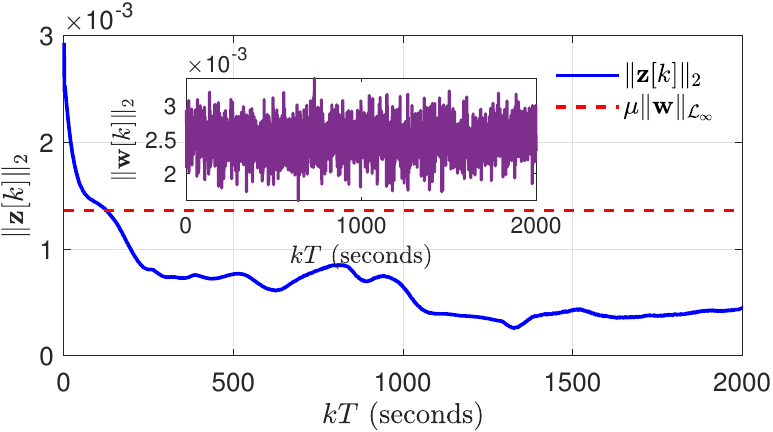}}{}\vspace{-0.1cm}\vspace{-0.2cm}
	\subfloat[\label{fig:scenario_7_performance_hwy2}]{\includegraphics[keepaspectratio=true,scale=0.61]{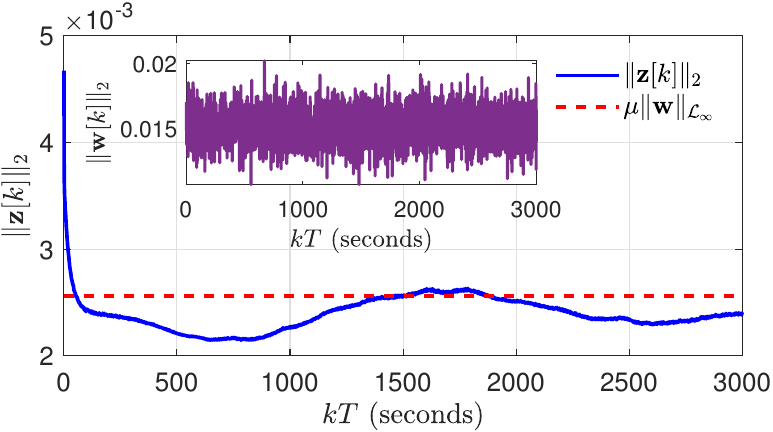}}{}\vspace{-0.1cm}\hspace{-0.0cm}\vspace{-0.02cm}
	\caption{Performance of estimation error with $40\%$ measurements for (a) Highway A (with determinant observability metric) and (b) Highway B (with trace observability metric).}
	\label{fig:scenario_7_8_perf}\vspace{-0.3cm}
\end{figure} 

%\begin{figure*}
%	\vspace{-0.0cm}
%	\centering 
%	\subfloat[]{\includegraphics[keepaspectratio=true,scale=0.5]{density_mainline_scenario_8}}{}\hspace{-0.08cm}
%	\subfloat[]{\includegraphics[keepaspectratio=true,scale=0.5]{density_onramps_scenario_8}}{}\hspace{-0.08cm}
%	\subfloat[]{\includegraphics[keepaspectratio=true,scale=0.5]{density_offramps_scenario_8}}{}\hspace{-0.08cm}\vspace{-0.00cm}
%	\caption{Comparison between real $\m x [k]$ and estimated $\hat{\m x} [k]$ traffic densities on (a) mainline segments, (b) on-ramps, and (c) off-ramps for highway A with $40\%$ measurements---sensor locations are obtained form solving \textbf{P1} with determinant observability metric.}
%	\label{fig:density_evolution_scenario_8}\vspace{-0.2cm}
%\end{figure*} 

\begin{figure*}
	\vspace{-0.40cm}
	\centering 
	\subfloat[]{\includegraphics[keepaspectratio=true,scale=0.5]{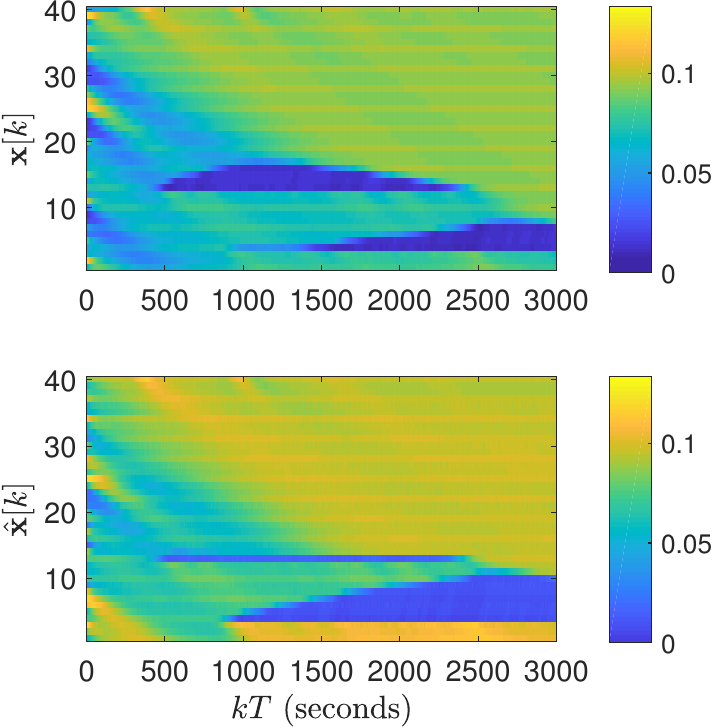}}{}\hspace{-0.08cm}
	\subfloat[]{\includegraphics[keepaspectratio=true,scale=0.5]{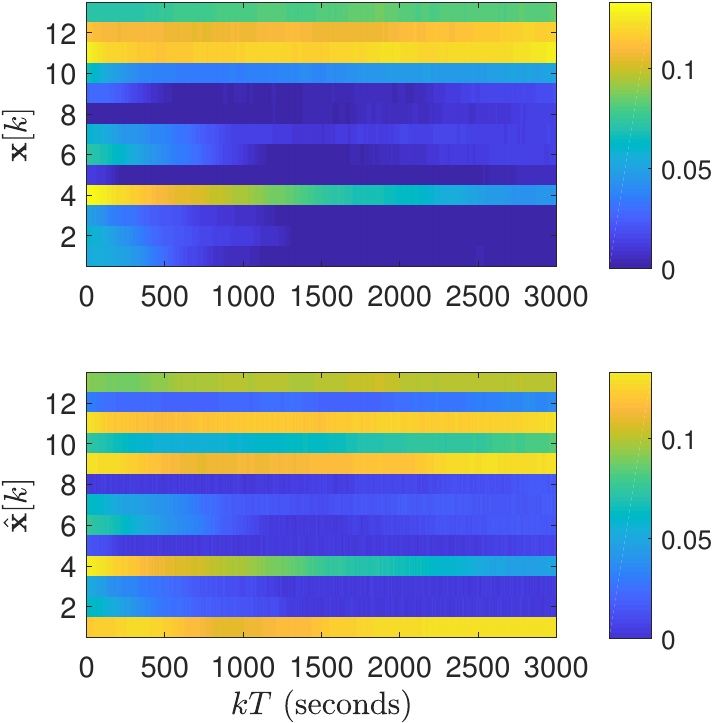}}{}\hspace{-0.08cm}
	\subfloat[]{\includegraphics[keepaspectratio=true,scale=0.5]{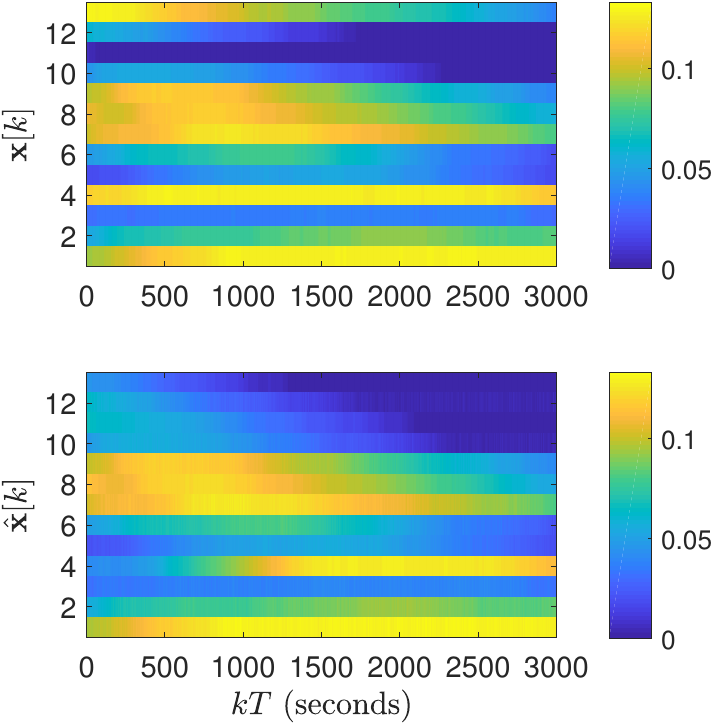}}{}\hspace{-0.08cm}\vspace{-0.00cm}
	\caption{Comparison between real $\m x [k]$ and estimated $\hat{\m x} [k]$ traffic densities on (a) mainline segments, (b) on-ramps, and (c) off-ramps for highway B with $40\%$ measurements---sensor locations are obtained form solving \textbf{P1} with trace observability metric.}
	\label{fig:density_evolution_scenario_7}\vspace{-0.1cm}
	\vspace{-0.2cm}
\end{figure*}   
 
\vspace{-0.2cm}

%\comment{Do you mention that the choice of the}
\subsection{Optimal, Randomized, and Uniform Sensor Placement }\label{ssec:optimal_vs_random}

%\comment{How does a uniform placement of sensors compare with our optimal and random? Like what happens if you select all even segments or all odd segments (that's what I mean by uniform? I'm curious. Let me know. Don't update the figures yet.) }
In the last part of numerical test, we analyze the results of traffic density estimation with optimal, randomized, and uniform sensor placement. We consider two highway networks: the first one is Highway A (detailed in Section \ref{ssec:obs_sensor_allocation}) while the other, referred to as \textit{Highway B}, consists a larger network with $40$ segments on the mainline, $13$ on-ramps, and $13$ off-ramps, giving $66$ highway segments in total with length of $9.9$ miles. The sensor placement and observer design for Highway A follow from two previous sections.
For Highway B, the sensor placement is acquired from solving \textbf{P1} with trace observability metric since it is much more efficient than the determinant, especially for larger network. The observation window is again set to be $N = 200$ with varying number of sensor allocations. We set the covariant matrix $R$ for Highway B to be $5\times10^{-3}\m I$ to simulate a higher measurement noise, while the final time is set to be $k_f = 3000$. The remaining parameters for observer design are chosen to remain the same as in Highway A. 

The first results of this numerical experiment are provided in Fig. \ref{fig:scenario_7_8_rmse}. 
This figure showcases the RMSE of traffic density estimation with optimal, randomized, and uniform sensor placement for both highways. 
For the uniform case, the sensors are firstly placed at odd locations and will begin to be placed at even locations only when all odd locations are used up.
Also note that, to compensate randomization on sensor placement, the simulations for this case (randomized sensor placement) are performed $10$ times for every $r$ and the results are averaged. For Highway A, the optimal sensor placements return significantly smaller RMSE than the randomized and uniform ones, regardless of observability metric being used. 
For Highway B, the RMSE for optimal sensor placement with trace is generally smaller than using randomized or uniform sensor placement.
These results suggest that using optimal sensor placement, regardless of
observability metric, gives better traffic density estimation than doing so with randomized sensor placement.  

The second results of this numerical experiment are illustrated in Fig. \ref{fig:scenario_7_8_perf}, where the corresponding error performance norm $\norm {\m z[k]}_2$ are compared with the estimation error performance bound $\mu\norm {\m w}_{\linf}$. Notice the larger value of error performance norm for Highway B due to higher measurement noise.
Specifically for both highways, it can be seen that the norm of error performance is converging to a value below the estimation error performance bound within the given simulation window. These findings are in accordance with the definition of $\linf$ stability. The resulting traffic density estimations for Highway B are depicted in Fig. \ref{fig:density_evolution_scenario_7}.
\vspace{-0.4cm}

\section{Conclusions, Paper Limitations, and Future Work}\label{sec:summary}

Given the thorough computational analysis in the previous section, the following observations are made, thereby answering the  posed research questions in Section \ref{sec:cases_tudy}:

\begin{enumerate}
	\item[--] \textit{A1:} Increasing observation window yields smaller initial state estimation errors while traffic network's quantifiable observability is directly proportional with observation window. This corroborates control-theoretic first principles. 
	\item[--]  \textit{A2:} Using trace as observability metric allows the sensor placement problem to be solved more efficiently than doing so with determinant metric for observability. This is critical if large-scale networks are considered.
	\item[--]  \textit{A3:} The resulting sensor placement obtained from solving \textbf{P1} is robust towards actual and presumed initial states.	This illustrates that the sensor placement to guarantee a wide range of initial operating conditions.
	\item[--]  \textit{A4:} When used to determine optimal sensor locations, the determinant is better than the trace in term of state estimation error and hence offers better network-wide observability. This poses tradeoffs between state estimation quality and computational tractability for the two metrics; see A2 above.  
	\item[--]  \textit{A5:} The optimal sensor placement outperforms randomized and uniform sensor placement in estimating traffic density. 
	\end{enumerate}
The approach in this paper has its own set of limitations. First, we consider a time-invariant traffic model where in reality, some parameters are in fact time-varying, which include critical density, split ratio, free-flow speed, and congestion wave speed. We also do not consider the capacity drop phenomenon or the stochasticity of traffic in the model presented in this paper.  Second, the sensor placement strategy does not consider the effect of measurement noise. To that end, future work will include solving the traffic sensor placement problem while considering a time-varying nonlinear traffic model incorporating capacity drop and stochasticity and  taking measurement noise into account, which is resulting in a robust sensor placement. Moreover, the submodularity properties of the determinant and trace---as well as other observability metrics---for sensor selection purpose will be investigated to assist in solving the placement problem for large  networks. Finally, we point out that the presented placement approach in the paper can be extended or utilized to other models in transportation systems beyond stretched highways, assuming that a nonlinear state-space representation of the dynamics is possible.

%\vspace{-0.144cm}

\section*{Acknowledgments}
We gratefully acknowledge the constructive comments from the reviewers---they have contributed positively to this work. We also acknowledge the financial support from the National Science Foundation through Grants 1636154, 1728629, 1739964, 1917164, 2152928, and 2152450.

\normalcolor

\bibliographystyle{IEEEtran}	\bibliography{arXiv_final_submitted}

\appendices

%\vspace{-0.245cm}

\section{State-Space Equation Parameters}\label{apdx:state_space_parameters}
In this section, we present the method for deriving the state-space equation \eqref{eq:state_space_gen}. For ease of reading, the time parameter is omitted from the notations for all the discrete time variables. Also, the state variables are written in terms of the traffic density variables that they represent. To obtain the equation \eqref{eq:state_space_gen}, we use the analytical expression for the non-linear $\min(\cdot,\cdot)$ function given as
\begin{equation}\label{eq:min_analytical}
    \min(a,b)=\dfrac{1}{2}\hspace{0.5mm}(\hspace{0.5mm}a\hspace{0.5mm}+\hspace{0.5mm}b\hspace{1mm}-\mid a\hspace{0.5mm}-\hspace{0.5mm}b\mid\hspace{0.5mm}).
\end{equation}
\noindent By applying \eqref{eq:min_analytical} to the non-linear state evolution equations \eqref{eq:non_linear_state_evolution} and simplifying, we can break the equations down to linear and non-linear parts that are lumped into the state-space equation parameters which are given below:

 \noindent $\m A \in \mathbb{R}^{n\times n}$, $n=N+N_I+N_O$. The elements of $\m A$ can be denoted by $a_{i,j}$ where $,i,j\in[1,n]$, and are given as follows
\begin{itemize}
    \item $i=1$
    \begin{align*}
        &a_{1,1} = 1-\dfrac{T}{4l}\left(v_f+w_c\right),  a_{1,2}=\dfrac{T}{4l}\left(w_c-\dfrac{\xi_{2}}{2}\right),\\ &a_{1,N+1}=\dfrac{T}{4l}v_f
    \end{align*}
    
    \item $i=N$
    \begin{align*}
    &a_{N,N-1} = \dfrac{T\bar{\beta}_{N-1}}{8l}v_f, a_{N,N} = 1-\dfrac{T}{4l}\left(v_f+w_c\right),\\ &a_{N,n}=-\dfrac{T\bar{\beta}_{N-1}}{4l\beta_{N-1}}w_c
    \end{align*}
    
    \item $i\in\Omega\setminus\Omega_I\cup\Omega_O,i\neq \{1,N\},j\in\hat{\Omega},k\in\check{\Omega}$
    \begin{align*}
    &a_{i,i-1}=\dfrac{T\bar{\beta}_{i-1}}{8l}v_f, a_{i,i}= 1-\dfrac{T}{4l}\left(v_f+w_c\right),\\ &a_{i,i+1}=\dfrac{T}{4l}\left(w_c-\dfrac{\xi_{i+1}}{2}\right), a_{i,N+j}=\dfrac{T}{4l}v_f,\\ &a_{i,N+N_I+k}=-\dfrac{T\bar{\beta}_{i-1}}{4l\beta_{i-1}}w_c
    \end{align*}
    
    \item $i\in\Omega_I,j\in\hat{\Omega}$
    \begin{align*}
    &a_{i,i-1}=\dfrac{T}{4l}v_f, a_{i,i}= 1-\dfrac{T}{4l}\left(v_f+w_c+\dfrac{\xi_i}{2}\right),\\ &a_{i,i+1}=\dfrac{T}{4l}w_c, a_{i,N+j}=\dfrac{T}{4l}v_f
    \end{align*}
    
    \item $i\in\Omega_O,j\in\check{\Omega}$
    \begin{align*}
    &a_{i,i-1}=\dfrac{T}{4l}v_f, a_{i,i}= 1-\dfrac{T}{4l}\left(\dfrac{1}{2}v_f+w_c\right),\\ &a_{i,i+1}=\dfrac{T}{4l\bar{\beta}_i}w_c, a_{i,N+N_I+j}=\dfrac{T}{4l\beta_i}w_c
    \end{align*}

    \item $i=N+\bar{i},\bar{i}\in\hat{\Omega},j\in{\Omega_I}$
    \begin{align*}
    &a_{i,i}=1-\dfrac{T}{2l}\left(v_f+\dfrac{w_c}{2}\right), a_{i,j}=\dfrac{T\xi_j}{4l}
    \end{align*}

    \item $i=N+N_I+\bar{i},\bar{i}\in\check{\Omega},j\in{\Omega_O}$
    \begin{align*}
    &a_{i,i}=1-\dfrac{T}{4l}\left(v_f+w_c\right), a_{i,j}=\dfrac{T\beta_j}{8l}v_f,\\
    &a_{i,j+1}=-\dfrac{T\beta_j}{4l\bar{\beta}_j}w_c
    \end{align*}

\end{itemize}

\noindent The remaining elements of $\m A$ that are not assigned a value above are equal to 0. 

 \noindent $\m f=\left[f_1\hspace{0.2cm}f_2\hspace{0.2cm}\dots \hspace{0.2cm}f_n\right]^T\in\mathbb{R}^{g}$ where $f_i,i\in[1,n]$ are column vectors of different but fixed lengths corresponding to the different cases of mainline sections and ramps, and $g$ is the sum of lengths of all the column vectors. The elements of each column vector $f_i$ can be denoted by $f_{j,i}$ where $j$ is in the range of the length of $f_i$, and are defined as follows:
\begin{itemize}
    \item $i=1$
    \begin{align*}
    &f_{1,1}=\dfrac{1}{w_cl}\mid f_{in}-\sigma_1\mid, f_{2,1}=\dfrac{1}{v_fw_cl}\mid \delta_1-\sigma_{2}\mid,\\ &f_{3,1}=\dfrac{1}{w_cl}\mid w_c\left(\rho_m-\rho_1\right)-v_f\rho_c\mid,\\ &f_{4,1}=\dfrac{1}{v_fl}\mid v_f\rho_{1}-v_f\rho_c\mid,\\ &f_{5,1}=\dfrac{1}{w_cl}\mid w_c\left(\rho_m-\rho_{2}\right)-v_f\rho_c\mid,\\ &f_{6,1}=\dfrac{1}{v_fw_cl}\mid v_f\hat{\rho}_{2}-\dfrac{\xi_{2}}{w_c}\min\left(w_c\left(\rho_m-\rho_{2}\right),v_f\rho_c\right)\mid\\
    &f_{7,1}=\dfrac{T}{4l}\left(2f_{in}+\left(\dfrac{\xi_{2}}{2w_c}-1\right)v_f\rho_c+\dfrac{\xi_{2}}{2}\rho_m\right)
    \end{align*}
    
    \item $i=N$
    \begin{align*}
    &f_{1,N}=\dfrac{1}{v_fw_cl}\mid \delta_{N-1}-\sigma_N\mid, f_{2,N}=\dfrac{1}{v_fl}\mid \delta_N-f_{out}\mid,\\ &f_{3,N}=\dfrac{1}{v_fw_cl}\mid \bar{\beta}_{N-1}\min\left(v_f\rho_{N-1},v_f\rho_c\right)\\&-\dfrac{\bar{\beta}_{N-1}}{\beta_{N-1}}\check{\sigma}_{N-1}\mid, f_{4,N}=\dfrac{1}{w_cl}\mid w_c\left(\rho_m-\rho_N\right)-v_f\rho_c\mid,\\ &f_{5,N}= \dfrac{1}{v_fl}\mid v_f\rho_{N}-v_f\rho_c\mid,\\ &f_{6,N}=\dfrac{1}{v_fl}\mid v_f\rho_{N-1}-v_f\rho_c\mid,\\ &f_{7,N}=\dfrac{T}{4l}\left(\dfrac{\bar{\beta}_{N-1}}{2}v_f\rho_c+\left(\dfrac{\bar{\beta}_{N-1}}{\beta_{N-1}}+1\right)w_c\rho_m\right)
    \end{align*}
    
    \item $i\in\Omega\setminus\Omega_I\cup\Omega_O,i\neq\{1,N\}$
    \begin{align*}
    &f_{1,i}=\dfrac{1}{v_fw_cl}\mid \delta_{i-1}-\sigma_i\mid, f_{2,i}=\dfrac{1}{v_fw_cl}\mid \delta_i-\sigma_{i+1}\mid,\\ &f_{3,i}=\dfrac{1}{v_fw_cl}\mid \bar{\beta}_{i-1}\min(v_f\rho_{i-1},v_f\rho_c)\\ &-\dfrac{\bar{\beta}_{i-1}}{\beta_{i-1}}w_c(\rho_m-\check{\rho}_{i-1})\mid,\\ &f_{4,i}=\dfrac{1}{w_cl}\mid w_c(\rho_m-\rho_i)-v_f\rho_c\mid,\\ &f_{5,i}=\dfrac{1}{v_fl}\mid v_f\rho_{i}-v_f\rho_c\mid,\\ &f_{6,i}=\dfrac{1}{w_cl}\mid w_c(\rho_m-\rho_{i+1})-v_f\rho_c\mid,\\ &f_{7,i}=\dfrac{1}{v_fw_cl}\mid v_f\hat{\rho}_{i+1}\\&-\dfrac{\xi_{i+1}}{w_c}\min(w_c(\rho_m-\rho_{i+1}),v_f\rho_c)\mid,\\ &f_{8,i}=\dfrac{1}{v_fl}\mid v_f\rho_{i-1}-v_f\rho_c\mid,\\ &f_{9,i}=\dfrac{T}{4l}\left(\left(\dfrac{\bar{\beta}_{i-1}}{2}+\dfrac{\xi_{i+1}}{2w_c}-1\right)v_f\rho_c+\left(\dfrac{\bar{\beta}_{i-1}}{\beta_{i-1}}\right.\right.\\&\left.\left.+\dfrac{\xi_{i+1}}{2w_c}\right)w_c\rho_m\right)
    \end{align*}
    
    \item $i\in\Omega_I$
    \begin{align*}
    &f_{1,i}=\dfrac{1}{v_fw_cl}\mid \delta_{i-1}-\sigma_i\mid, f_{2,i}=\dfrac{1}{v_fw_cl}\mid \delta_i-\sigma_{i+1}\mid,\\ &f_{3,i}=\dfrac{1}{v_fl}\mid v_f\rho_{i-1}-v_f\rho_c\mid,\\ &f_{4,i}=\dfrac{1}{w_cl}\mid w_c(\rho_m-\rho_i)-v_f\rho_c\mid,\\ &f_{5,i}=\dfrac{1}{v_fl}\mid v_f\rho_i-v_f\rho_c\mid,\\ &f_{6,i}=\dfrac{1}{w_cl}\mid w_c(\rho_m-\rho_{i+1})-v_f\rho_c\mid,\\ &f_{7,i}=\dfrac{1}{v_fw_cl}\mid v_f\hat{\rho}_i-\dfrac{\xi_{i}}{w_c}\min(w_c(\rho_m-\rho_{i}),v_f\rho_c)\mid,\\ &f_{8,i}=\dfrac{1}{w_cl}\mid w_c(\rho_m-\rho_{i})-v_f\rho_c\mid\\
    &f_{9,i}=\dfrac{T\xi_i}{8l}\left(\rho_m+\dfrac{v_f\rho_c}{w_c}\right)
    \end{align*}
    
    \item $i\in\Omega_O$
    \begin{align*}
    &f_{1,i}=\dfrac{1}{v_fw_cl}\mid \delta_{i-1}-\sigma_i\mid, f_{2,i}=\dfrac{1}{v_fw_cl}\mid \delta_i-\sigma_{i+1}\mid,\\  &f_{3,i}=\dfrac{1}{v_fl}\mid v_f\rho_{i-1}-v_f\rho_c\mid,\\ &f_{4,i}=\dfrac{1}{w_cl}\mid w_c(\rho_m-\rho_i)-v_f\rho_c\mid,\\ &f_{5,i}=\dfrac{1}{v_fw_cl}\mid \min(v_f\rho_i,v_f\rho_c)-\dfrac{1}{\beta_i}w_c(\rho_m-\check{\rho}_i) \mid,\\ &f_{6,i}=\dfrac{1}{w_cl}\mid w_c(\rho_m-\rho_{i+1})-v_f\rho_c \mid,\\ &f_{7,i}=\dfrac{1}{v_fl}\mid v_f\rho_i-v_f\rho_c\mid,\\ &f_{8,i}=\dfrac{T}{4l}\left(\left(\dfrac{3}{2}-\dfrac{1}{\bar{\beta}_i}\right)v_f\rho_c+\left(1-\dfrac{1}{\beta_i}-\dfrac{1}{\bar{\beta}_i}\right)w_c\rho_m\right)
    \end{align*}
    
    \item $i=N+\bar{i},\bar{i}\in\hat{\Omega},j\in{\Omega_I}$
    \begin{align*}
    &f_{1,i}=\dfrac{1}{w_cl}\mid \hat{\sigma}_j-\hat{f}_j \mid,\\ &f_{2,i}=\dfrac{1}{v_fw_cl}\mid v_f\hat{\rho}_{j}-\dfrac{\xi_{j}}{w_c}\min(w_c(\rho_m-\rho_{j}),v_f\rho_c)\mid,\\ &f_{3,i}=\dfrac{1}{w_cl}\mid w_c(\rho_m-\hat{\rho}_j)-v_f\rho_c\mid,\\ &f_{4,i}=\dfrac{1}{w_cl}\mid w_c(\rho_m-\rho_{j})-v_f\rho_c\mid,\\ &f_{5,i}=\dfrac{T}{2l}\left(\left(\dfrac{1}{2}-\dfrac{\xi_j}{2w_c}\right)v_f\rho_c+\left(\dfrac{1}{2}-\dfrac{\xi_j}{2w_c}\right)w_c\rho_m\right)
    \end{align*}
    
    \item $i=N+N_I+\bar{i},\bar{i}\in\check{\Omega},j\in{\Omega_O}$
    \begin{align*}
    &f_{1,i}=\dfrac{1}{v_fw_cl}\mid \delta_j - \sigma_{j+1}\mid, f_{2,i}=\dfrac{1}{v_fl}\mid \check{\delta}_j-\check{f}_j\mid,\\ &f_{3,i}=\dfrac{1}{v_fw_cl}\mid \min(v_f\rho_j,v_f\rho_c)-\dfrac{1}{\beta_j}\check{\sigma}_j\mid,\\ &f_{4,i}=\dfrac{1}{w_cl}\mid w_c(\rho_m-\rho_{j+1})-v_f\rho_c\mid,\\ &f_{5,i}=\dfrac{1}{v_fl}\mid v_f\check{\rho}_j-v_f\rho_c\mid,f_{6,i}=\dfrac{1}{v_fl}\mid v_f\rho_j-v_f\rho_c\mid,\\ &f_{7,i}=\dfrac{T}{2l}\left(\left(\dfrac{\beta_j}{4}+\dfrac{\beta_j}{2\bar{\beta}_j}-\dfrac{1}{2}\right)v_f\rho_c+\left(\dfrac{1}{2}+\dfrac{\beta_j}{2\bar{\beta}_j}\right)w_c\rho_m\right)
    \end{align*} 
\end{itemize}
\vspace{0.1cm}

\noindent $\m G\in\mathbb{R}^{n\times g}$ is of the form\\
\\
$
\begin{bmatrix}G_1 & 0  & \cdots & \cdots  & 0\\
0 & G_2 & \ddots &  & \vdots\\
\vdots & \ddots &\ddots & \ddots &\vdots \\
\vdots & & \ddots& G_{n-1} & 0\\ 0 & \cdots& \cdots& 0&G_{n}&
\end{bmatrix},$
\vspace{0.1cm}
where $G_i, i\in[1,n]$ are row vectors with the same number of columns as the number of rows in $f_i$. The elements of each row vector $G_i$ can be denoted by $G_{i,j}$, where $j$ is in the range of the length of $G_i$, and are defined as follows:
\begin{itemize}
    \item $i=1$
    \begin{multline*}
     G_i = \dfrac{Tv_fw_c}{4}\left[\begin{matrix} -\dfrac{2}{v_f} & 2 & -\dfrac{1}{v_f} & \dfrac{1}{w_c} \end{matrix} \right. \\
     \left.\begin{matrix} \left(\dfrac{1}{v_f}-\dfrac{\xi_2}{2v_fw_c}\right) & -1 & \dfrac{4}{Tv_fw_c}\end{matrix}\right]
    \end{multline*}
    % \begin{align*}
    %     &G_{1,1}-\dfrac{Tw_c}{2}, G_{1,2}=\dfrac{Tv_fw_c}{2}, G_{1,3}-\dfrac{Tw_c}{4}, G_{1,4}\dfrac{Tv_f}{4},\\ &G_{1,5}(\dfrac{Tw_c}{4}-\dfrac{T\xi_{2}}{8}),G_{1,6}-\dfrac{Tv_fw_c}{4},G_{1,7}=1
    % \end{align*}
    
    \item $i=N$
    \begin{multline*}
     G_i = \dfrac{Tv_fw_c}{4}\left[\begin{matrix} -2 & \dfrac{2}{w_c} & -1 & -\dfrac{1}{v_f} & \dfrac{1}{w_c}\end{matrix} \right. \\
     \left.\begin{matrix}  & -\dfrac{\bar{\beta}_{N-1}}{2w_c} & \dfrac{4}{Tv_fw_c}\end{matrix}\right]
    \end{multline*}
    % \begin{align*}
    % &G_{N,1}=-\dfrac{Tv_fw_c}{2}, G_{N,2}=\dfrac{Tv_f}{2}, G_{N,3}=-\dfrac{Tv_fw_c}{4},\\ &G_{N,4}=-\dfrac{Tw_c}{4}, G_{N,5}=\dfrac{Tv_f}{4},G_{N,6}=-\dfrac{Tv_f\bar{\beta}_{N-1}}{8},\\&G_{N,7}=1
    % \end{align*}

    \item $i\in\Omega\setminus\Omega_I\cup\Omega_O,i\neq\{1,N\}$\\
    % \begin{align*}
    \begin{multline*}
    G_i = \dfrac{Tv_fw_c}{4}\left[\begin{matrix} -2 & 2 & -1 & -\dfrac{1}{v_f} & \dfrac{1}{w_c}\end{matrix}\right. \\
     \left.\begin{matrix} \left(\dfrac{1}{v_f}-\dfrac{\xi_{i+1}}{2v_fw_c}\right) & -1 & -\dfrac{\bar{\beta}_{i-1}}{2w_c} & \dfrac{4}{Tv_fw_c}\end{matrix}\right]
    \end{multline*}
    % &G_{i,1}=-\dfrac{Tv_fw_c}{2}, G_{i,2}=\dfrac{Tv_fw_c}{2}, G_{i,3}-\dfrac{Tv_fw_c}{4},\\ &G_{i,4}=-\dfrac{Tw_c}{4}, G_{i,5}=\dfrac{Tv_f}{4}, G_{i,6}=\left(\dfrac{Tw_c}{4}-\dfrac{T\xi_{i+1}}{8}\right),\\ &G_{i,7}=-\dfrac{Tv_fw_c}{4},G_{i,8}=-\dfrac{Tv_f\bar{\beta}_{i-1}}{8}, G_{i,9}=1
    % \end{align*}
    
    \item $i\in\Omega_I$\\
    % \begin{align*}
    \begin{multline*}
     G_i = \dfrac{Tv_fw_c}{4}\left[\begin{matrix} -2 & 2 & -\dfrac{1}{w_c} & -\dfrac{1}{v_f} & \dfrac{1}{w_c} &\dfrac{1}{v_f} \end{matrix} \right. \\
     \left.\begin{matrix}  -1 & -\dfrac{\xi_i}{2v_fw_c} & \dfrac{4}{Tv_fw_c}\end{matrix}\right]
    \end{multline*}
    % &G_{i,1}=-\dfrac{Tv_fw_c}{2}, G_{i,2}=\dfrac{Tv_fw_c}{2}, G_{i,3}=-\dfrac{Tv_f}{4},\\ &G_{i,4}=-\dfrac{Tw_c}{4}, G_{i,5}=\dfrac{Tv_f}{4}, G_{i,6}=G_{i,}=\dfrac{Tw_c}{4},\\ &G_{i,7}=-\dfrac{Tv_fw_c}{4}, G_{i,8}=-\dfrac{T\xi_i}{8}, G_{i,9}=1
    % \end{align*}

    \item $i\in\Omega_O$
    \begin{multline*}
     G_i = \dfrac{Tv_fw_c}{4}\left[\begin{matrix} -2 & \dfrac{2}{\bar{\beta}_i} & -\dfrac{1}{w_c} & -\dfrac{1}{v_f} & 1 \end{matrix} \right. \\
     \left.\begin{matrix} \dfrac{1}{v_f\bar{\beta}_i} & \dfrac{1}{2w_c} & \dfrac{4}{Tv_fw_c}\end{matrix}\right]
    \end{multline*}
    % \begin{align*}
    % &G_{i,1}=-\dfrac{Tv_fw_c}{2}, G_{i,2}=\dfrac{Tv_fw_c}{2\bar{\beta}_i}, G_{i,3}=-\dfrac{Tv_f}{4},\\ &G_{i,4}=-\dfrac{Tw_c}{4}, G_{i,5}=\dfrac{Tv_fw_c}{4}, G_{i,6}=\dfrac{Tw_c}{4\bar{\beta}_i},\\&G_{i,7}=\dfrac{Tv_f}{8}, G_{i,8} = 1
    % \end{align*}

    \item $i=N+\bar{i},\bar{i}\in\hat{\Omega},j\in{\Omega_I}$
    \begin{align*}
     G_i = \dfrac{Tw_c}{4}\left[\begin{matrix} -2 & 2v_f & -1 & \dfrac{\xi_j}{w_c} & \dfrac{4}{Tw_c} \end{matrix} \right]
    \end{align*}
    % \begin{align*}
    % &G_{i,1}=-\dfrac{Tw_c}{2}, G_{i,2}=\dfrac{Tv_fw_c}{2}, G_{i,3}=-\dfrac{Tw_c}{4},\\ &G_{i,4}=\dfrac{T\xi_j}{4},G_{i,5}=1
    % \end{align*}

    \item $i=N+N_I+\bar{i},\bar{i}\in\check{\Omega},j\in{\Omega_O}$
    \begin{multline*}
     G_i = \dfrac{Tv_f}{4}\left[\begin{matrix} -\dfrac{2w_c\beta_j}{\bar{\beta}_j} & 2 & -w_c\beta_j & -\dfrac{w_c\beta_j}{v_f\bar{\beta}_j} & 1 \end{matrix} \right. \\
     \left.\begin{matrix} -\dfrac{\beta_j}{2} & \dfrac{T}{4v_f}\end{matrix}\right]
    \end{multline*}
    % \begin{align*}
    % &G_{i,1}=-\dfrac{Tv_fw_c\beta_j}{2\bar{\beta}_j}, G_{i,2}=\dfrac{Tv_f}{2},G_{i,3}=-\dfrac{Tv_fw_c\beta_j}{4},\\ &G_{i,4}=-\dfrac{Tw_c\beta_j}{4\bar{\beta}_j}, G_{i,5}=\dfrac{Tv_f}{4}, G_{i,6}=-\dfrac{Tv_f\beta_j}{8},\\&G_{i,7}=1
    % \end{align*} 
\end{itemize}
$\m {B_u}\in\mathbb{R}^{n\times(2+N_I+N_O)}$. The elements of $\m{B_u}$ can be denoted by $\m B_{i,j}$ and are given by the following function\\
$\m {B}_{i,j}=
\begin{cases}
\dfrac{T}{2l}f_{in},&\mathrm{if} \:i=j=1\\
-\dfrac{T}{2l}f_{out},&\mathrm{if} \:i=N,j=2\\
\dfrac{T}{2l}\hat{f}_k,&\mathrm{if} \:i=N+\bar{i},j=2+\bar{i},\bar{i}\in\hat{\Omega},k\in\Omega_I\\
-\dfrac{T}{2l}\check{f}_k,&\mathrm{if} \:i=N+N_I+\bar{i},j=2+N_I+\bar{i},\\
&\bar{i}\in\check{\Omega},k\in\Omega_I\\
0, &\mathrm{otherwise}.
\end{cases}
$
%\vspace{-0.4cm}
\section{Calculation of Jacobian Matrix}
\label{apdx:calc_jacobian}
To calculate the Jacobian matrix for the non-linear vector valued function $\m f(\m x,\m u)$ containing $\mid\cdot\mid$ and $\min(\cdot,\cdot)$ functions, we first replace all the $\min(\cdot,\cdot)$ functions by their analytical form \eqref{eq:min_analytical} and then calculate the derivative using the definition of the derivative for $\mid\cdot\mid$ function given as 
\begin{align}\label{e:differential_absolute}
    \dfrac{\partial \mid x\mid}{\partial x}=
    \dfrac{x}{\mid x \mid},x\neq 0
\end{align}
At $x=0$, we simply take the derivative to be 0. The elements of the Jacobian matrix can be calculated using \eqref{e:differential_absolute} along with the chain rule of differentiation. For example, the element $f_{2,i}$ in $f_i,i\in \Omega_I,j=N+\bar{j},\bar{j}\in\hat{\Omega}$ can be differentiated with respect to the state variable $x_i$, for $\delta_{i-1}\neq\sigma_{i}$, as follows:
\begin{align*}
\dfrac{\partial f_{1,i}}{\partial x_i} &= \dfrac{1}{v_fw_cl}\cdot \dfrac{\partial \mid \delta_{i-1}-\sigma_{i}\mid}{\partial x_i}\\
% &=\frac{1}{v_fw_cl}\frac{\delta_{i-1}-\sigma_{i}}{\mid \delta_{i-1}-\sigma_{i}\mid}\frac{\partial (\delta_{i-1}-\sigma_{i})}{\partial x_i}\\
&=-\dfrac{1}{v_fw_cl}\cdot\dfrac{\delta_{i-1}-\sigma_{i}}{\mid \delta_{i-1}-\sigma_{i}\mid}\cdot\dfrac{\partial \sigma_{i}}{\partial x_i},
\end{align*}
% where, in terms of the state vector $\m x$
where
\begin{align*}
\dfrac{\partial \sigma_{i}}{\partial x_i}
&=\dfrac{\partial \min(w_c(\rho_m-x_i),v_f\rho_c)}{\partial x_i}-\dfrac{\partial r_i}{\partial x_i}.
\end{align*}
% where
Here, for $w_c(\rho_m-x_i)\neq v_f\rho_c$,
\begin{align*}
&\dfrac{\partial \min(w_c(\rho_m-x_i),v_f\rho_c)}{\partial x_i}\\
% &=\frac{1}{2}\frac{\partial (w_c(\rho_m-x_i)+v_f\rho_c-\mid w_c(\rho_m-x_i)-v_f\rho_c \mid)}{\partial x_i}\\
% &=-\frac{1}{2}w_c+\frac{1}{2}\frac{\partial(-\mid w_c(\rho_m-x_i)-v_f\rho_c \mid)}{\partial x_i}\\
% &=-\frac{1}{2}w_c+\frac{1}{2}w_c\frac{w_c(\rho_m-x_i)-v_f\rho_c}{\mid w_c(\rho_m-x_i)-v_f\rho_c \mid}\\
&=-\dfrac{w_c}{2}\left(1-\dfrac{w_c(\rho_m-x_i)-v_f\rho_c}{\mid w_c(\rho_m-x_i)-v_f\rho_c \mid}\right),
\end{align*}
and for $v_fx_j\neq\dfrac{\xi_{i}}{w_c}\min(w_c(\rho_m-x_i),v_f\rho_c)$,
\begin{align*}
\dfrac{\partial r_i}{\partial x_i}&=
% \dfrac{\partial \min\left(v_fx_j,\dfrac{\xi_{i}}{w_c}\min(w_c(\rho_m-x_i),v_f\rho_c)\right)}{\partial x_i}\\
% &=\frac{1}{2}\left(\frac{\partial (v_fx_j+\frac{\xi_{i}}{w_c}\min(w_c(\rho_m-x_i),v_f\rho_c))}{\partial x_i}\\
% &-\frac{\partial (\mid v_fx_j-\frac{\xi_{i}}{w_c}\min(w_c(\rho_m-x_i),v_f\rho_c)\mid)}{\partial x_i}\right)\\
% &=\frac{1}{2}\left(\frac{\xi_i}{w_c}\frac{\partial \min(w_c(\rho_m-x_i),v_f\rho_c)}{\partial x_i}\\
% &+\frac{v_fx_j-\frac{\xi_{i}}{w_c}\min(w_c(\rho_m-x_i),v_f\rho_c)}{\mid v_fx_j-\frac{\xi_{i}}{w_c}\min(w_c(\rho_m-x_i),v_f\rho_c)\mid}\cdot\\
% &\frac{\xi_i}{w_c}\frac{\partial\min(w_c(\rho_m-x_i),v_f\rho_c)}{\partial x_i}\right)\\
% &=\frac{\xi_i}{2w_c}\frac{\partial\min(w_c(\rho_m-x_i),v_f\rho_c)}{\partial x_i}\left(1\\
% &+\frac{v_fx_j-\frac{\xi_{i}}{w_c}\min(w_c(\rho_m-x_i),v_f\rho_c)}{\mid v_fx_j-\frac{\xi_{i}}{w_c}\min(w_c(\rho_m-x_i),v_f\rho_c)\mid}\right)\\
-\dfrac{\xi_i}{2w_c}\cdot\dfrac{w_c}{2}\left(1-\dfrac{w_c(\rho_m-x_i)-v_f\rho_c}{\mid w_c(\rho_m-x_i)-v_f\rho_c \mid}\right)\cdot\\
&\left(1+\dfrac{v_fx_j-\dfrac{\xi_{i}}{w_c}\min(w_c(\rho_m-x_i),v_f\rho_c)}{\mid v_fx_j-\dfrac{\xi_{i}}{w_c}\min(w_c(\rho_m-x_i),v_f\rho_c)\mid}\right).
\end{align*}

\noindent Therefore, the derivative can be written as
\begin{align*}
    &\dfrac{\partial f_{1,i}}{\partial x_i}=\dfrac{1}{2v_fl}\cdot\dfrac{\delta_{i-1}-\sigma_{i}}{\mid \delta_{i-1}-\sigma_{i}\mid}\left(1-\dfrac{w_c(\rho_m-x_i)-v_f\rho_c}{\mid w_c(\rho_m-x_i)-v_f\rho_c \mid}\right)\cdot\\
    &\left(1-\dfrac{\xi_i}{2w_c}\left(1+\dfrac{v_fx_j-\dfrac{\xi_{i}}{w_c}\min(w_c(\rho_m-x_i),v_f\rho_c)}{\mid v_fx_j-\dfrac{\xi_{i}}{w_c}\min(w_c(\rho_m-x_i),v_f\rho_c)\mid}\right)\right).
\end{align*}
\noindent Note that the above expression is valid when $\delta_{i-1}\neq\sigma_{i}$, $w_c(\rho_m-x_i)\neq v_f\rho_c$ and $v_fx_j\neq\dfrac{\xi_{i}}{w_c}\min(w_c(\rho_m-x_i),v_f\rho_c)$. If, on the other hand,  $w_c(\rho_m-x_i)= v_f\rho_c$, then
\begin{align*}
    &\dfrac{\partial f_{1,i}}{\partial x_i}=\dfrac{1}{2v_fl}\cdot\dfrac{\delta_{i-1}-\sigma_{i}}{\mid \delta_{i-1}-\sigma_{i}\mid}\cdot\\
    &\left(1-\dfrac{\xi_i}{2w_c}\left(1+\dfrac{v_fx_j-\dfrac{\xi_{i}}{w_c}\min(w_c(\rho_m-x_i),v_f\rho_c)}{\mid v_fx_j-\dfrac{\xi_{i}}{w_c}\min(w_c(\rho_m-x_i),v_f\rho_c)\mid}\right)\right),
\end{align*}
\\
if $v_fx_j=\dfrac{\xi_{i}}{w_c}\min(w_c(\rho_m-x_i),v_f\rho_c)$, then\\
\begin{align*}
    \dfrac{\partial f_{1,i}}{\partial x_i}=&\dfrac{1}{2v_fl}\cdot\dfrac{\delta_{i-1}-\sigma_{i}}{\mid \delta_{i-1}-\sigma_{i}\mid}\left(1-\dfrac{\xi_i}{2w_c}\right)\cdot\\
    &\left(1-\dfrac{w_c(\rho_m-x_i)-v_f\rho_c}{\mid w_c(\rho_m-x_i)-v_f\rho_c \mid}\right).
\end{align*}

\noindent if both $w_c(\rho_m-x_i)= v_f\rho_c$, and $v_fx_j=\dfrac{\xi_{i}}{w_c}\min(w_c(\rho_m-x_i),v_f\rho_c)$, then
\begin{align*}
    \dfrac{\partial f_{1,i}}{\partial x_i}=&\dfrac{1}{2v_fl}\cdot\dfrac{\delta_{i-1}-\sigma_{i}}{\mid \delta_{i-1}-\sigma_{i}\mid}\left(1-\dfrac{\xi_i}{2w_c}\right).
\end{align*}

\noindent and if $\delta_{i-1}=\sigma_{i}$, then
%\begin{align*}
$\dfrac{\partial f_{1,i}}{\partial x_i}=0.$
%\end{align*}
Similarly, other elements of the Jacobian matrix can also be derived. 

\section{Calculation of Lipschitz constant}
\label{apdx:calc_lipschitz_constant}
In this section we present the methodology for calculating a Lipschitz constant $\gamma$ that satisfies
\begin{align*}
    \lvert\lvert \m f(\m x,u)-\m f(\hat{\m x},u)\rvert\rvert_2\le\gamma\lvert\lvert \m x-\m {\hat{x}}\rvert\rvert_2
\end{align*}

To calculate $\gamma$, we use the following inequalities related to the $\mid\cdot\mid$ function:
\begin{subequations}
\begin{align}
    &\mid\mid a\mid-\mid b\mid\mid\hspace{2mm}\le\hspace{2mm}\mid a-b\mid\label{eq:ineq1}\\
    &\mid a \pm b\mid\hspace{2mm}\le\hspace{2mm}\mid a\mid+\mid b\mid\label{eq:ineq2}\\
    &\mid a_1\mid+\mid a_2\mid\hspace{2mm}\le \hspace{2mm}\sqrt{2}(a_1^2+a_2^2)^{1/2}=\sqrt{2}\lvert\lvert \m a\rvert\rvert_2\label{eq:ineq3}
\end{align}
\end{subequations}
By using one or more of these inequalities we can calculate the Lipschitz constant for each component of the vector function $\m f(\cdot)$, and combine them to obtain the final $\gamma$. 

For example, let us calculate the the Lipschitz constant for the expression $f_{1,i},i\in\Omega_I$, 
% given as:
% \begin{align*}
% &f_{1,i}=\frac{1}{v_fw_cl}\mid \delta_{i-1}-\sigma_i\mid\\ 
% &=\frac{1}{v_fw_cl}\mid \min(v_f\rho_{i-1},v_f\rho_c) - \min(w_c(\rho_m-\rho_i),v_f\rho_c)-r_i\mid\\
% &=\frac{1}{v_fw_cl}\mid \min(v_f\rho_{i-1},v_f\rho_c) - \min(w_c(\rho_m-\rho_i),v_f\rho_c)\\&-\min\big(v_f\hat{\rho}_{i},\frac{\xi_{i}}{w_c}\min(w_c(\rho_m-\rho_{i}),v_f\rho_c)\big)\mid
% \end{align*}
% To avoid any confusion between the actual and estimated state vectors, especially for the on-ramp densities, we replace the $\rho_{i-1}$,$\rho_i$ and $\hat{\rho}_i$ variables in the above equation with their state vector equivalents.
% This means that the above equation, for
% $i\in\Omega_I$ and 
which, for $j=N+\bar{j},\bar{j}\in\hat{\Omega}$, is given as
\begin{align*}
    &f_{1,i}=\frac{1}{v_fw_cl}\mid \min(v_fx_{i-1},v_f\rho_c) - \min(w_c(\rho_m-x_i),v_f\rho_c)\\&-\min\big(v_fx_j,\frac{\xi_{i}}{w_c}\min(w_c(\rho_m-x_{i}),v_f\rho_c)\big)\mid.
\end{align*}

% Using \eqref{eq:min_analytical} and the fact that $v_f\rho_c=w_c(\rho_m-\rho_c)$, and simplifying, we get
% \begin{align*}
%     &f_{1,i}=\frac{1}{2v_fw_cl}\mid v_fx_{i-1}+ \Big(1-\frac{\xi_i}{2w_c}\Big)w_cx_i+ v_fx_j\\ &+w_c\rho_m\Big(\frac{\xi_i}{2w_c}-1\Big)+\frac{\xi_iv_f\rho_c}{2w_c}-\mid v_fx_{i-1}-v_f\rho_c\mid-\mid v_fx_{j}\\&-\frac{\xi_{i}}{w_c}\min(w_c(\rho_m-x_{i}),v_f\rho_c)\mid+\Big(1-\frac{\xi_i}{2w_c}\Big)\mid w_cx_i\\&-w_c\rho_c\mid \mid
% \end{align*}
% Applying \eqref{eq:ineq1}, we can write the inequality for $\mid f_{1,i}-\hat{f}_{1,i}\mid$ as follows
% \begin{align*}
%     &\mid f_{1,i}-\hat{f}_{1,i}\mid\le
%     \frac{1}{2v_fw_cl}\mid (v_fx_{i-1}-v_f\hat{x}_{i-1})\\&+ \Big(1-\frac{\xi_i}{2w_c}\Big)(w_cx_i-w_c\hat{x}_i)+ (v_fx_j-v_f\hat{x}_j)\\ &-(\mid v_fx_{i-1}-v_f\rho_c\mid-\mid v_f\hat{x}_{i-1}-v_f\rho_c\mid)\\&-(\mid v_fx_{j}-\frac{\xi_{i}}{w_c}\min(w_c(\rho_m-x_{i}),v_f\rho_c)\mid\\&-\mid v_f\hat{x}_{j}-\frac{\xi_{i}}{w_c}\min(w_c(\rho_m-\hat{x}_{i}),v_f\rho_c)\mid)\\&+\Big(1-\frac{\xi_i}{2w_c}\Big)(\mid w_cx_i-w_c\rho_c\mid-\mid w_c\hat{x}_i-w_c\rho_c\mid) \mid
% \end{align*}
% Again applying \eqref{eq:ineq1}, we get

Using \eqref{eq:min_analytical}, \eqref{eq:ineq1}, and \eqref{eq:ineq2}, and simplifying, we can write the following inequality for $\mid f_{1,i}-\hat{f}_{1,i}\mid$:
% \begin{align*}
%     &\mid f_{1,i}-\hat{f}_{1,i}\mid\le
%     \frac{1}{2v_fw_cl}\mid (v_fx_{i-1}-v_f\hat{x}_{i-1})\\&+ \Big(1-\frac{\xi_i}{2w_c}\Big)(w_cx_i-w_c\hat{x}_i)+ (v_fx_j-v_f\hat{x}_j)\\ &-\mid v_fx_{i-1}-v_f\hat{x}_{i-1}\mid-\mid v_fx_{j}\\&-\frac{\xi_{i}}{w_c}\min(w_c(\rho_m-x_{i}),v_f\rho_c)- v_f\hat{x}_{j}\\&+\frac{\xi_{i}}{w_c}\min(w_c(\rho_m-\hat{x}_{i}),v_f\rho_c)\mid\\&+\Big(1-\frac{\xi_i}{2w_c}\Big)\mid w_cx_i-w_c\hat{x}_i\mid \mid
% \end{align*}

% Applying \eqref{eq:ineq2}, we can simplify this to
% \begin{align*}
%     &\mid f_{1,i}-\hat{f}_{1,i}\mid\le
%     \frac{1}{2v_fw_cl}(v_f\mid x_{i-1}-\hat{x}_{i-1}\mid\\&+ \Big(1-\frac{\xi_i}{2w_c}\Big)w_c\mid x_i-\hat{x}_i\mid+ v_f\mid x_j-\hat{x}_j\mid\\ &+v_f\mid x_{i-1}-\hat{x}_{i-1}\mid+v_f\mid x_{j}-\hat{x}_j\mid\\&+\mid\frac{\xi_{i}}{w_c}\min(w_c(\rho_m-x_{i}),v_f\rho_c)\\&-\frac{\xi_{i}}{w_c}\min(w_c(\rho_m-\hat{x}_{i}),v_f\rho_c)\mid\\&+\Big(1-\frac{\xi_i}{2w_c}\Big)w_c\mid x_i-\hat{x}_i\mid)
% \end{align*}
% Again using \eqref{eq:min_analytical} and inequalities as above and simplifying, we get
\begin{align*}
    &\mid f_{1,i}-\hat{f}_{1,i}\mid\le
    \frac{1}{2v_fw_cl}(2v_f\mid x_{i-1}-\hat{x}_{i-1}\mid\\&+ 2\Big(1-\frac{\xi_i}{2w_c}\Big)w_c\mid x_i-\hat{x}_i\mid+ 2v_f\mid x_j-\hat{x}_j\mid\\&+\xi_i\mid x_i-\hat{x}_i\mid).
\end{align*}

By applying \eqref{eq:ineq3} to combine the terms with the same coefficients and then adding up all the coefficients, we can write the above inequality as
\begin{align*}
    % &\mid f_{1,i}-\hat{f}_{1,i}\mid\le
    % \frac{1}{l}\Big(\frac{\sqrt{2}}{w_c}+\frac{1}{v_f}\Big(1-\frac{\xi_i}{2w_c}\Big)+\frac{\xi_i}{2v_fw_c}\Big)\lvert\lvert \m x-\m {\hat{x}}\rvert\rvert_2\\
    % &\implies 
    \mid f_{1,i}-\hat{f}_{1,i}\mid\le
    \frac{1}{l}\Big(\frac{\sqrt{2}}{w_c}+\frac{1}{v_f}\Big)\lvert\lvert \m x-\m {\hat{x}}\rvert\rvert_2.
\end{align*}

Let us call the coefficient of $\lvert\lvert \m x-\m {\hat{x}}\rvert\rvert_2$ in the above equation as $\gamma_{1,i}$. We can obtain the $\gamma_i$ corresponding to $f_i$ by adding these coefficient values for all the components of the vector $f_i$.
% For some expressions like $f_{5,i},i\in\Omega_I$, $\gamma_i$ is much more obvious where we obtain the following inequality
% \begin{align*}
%     &\mid f_{5,i}-\hat{f}_{5,i}\mid\le\frac{1}{l}\mid x_i-\hat{x}_i\mid
% \end{align*}

The values of $\gamma_i$ for different cases of mainline sections and ramps are given below:
\begin{itemize}
    \item $i=1$
    \begin{align*}
        \frac{1}{l}\Big(\frac{1+\sqrt{2}}{w_c} + \frac{1}{v_f} + \frac{\xi_2}{v_fw_c}+4\Big)
    \end{align*}
    
    \item $i=N$
    \begin{align*}
        \frac{1}{l}\Big(\frac{2\bar{\beta}_{N-1}}{w_c} + \Big(1+2\frac{\bar{\beta}_{N-1}}{\beta_{N-1}}\Big)\frac{1}{v_f}+4\Big)
    \end{align*}
    
    \item $i\in\Omega\setminus\Omega_I\cup\Omega_O,i\neq\{1,N\}$
    \begin{align*}
        \frac{1}{l}\Big(\frac{1+\sqrt{2}+2\bar{\beta}_{i-1}}{w_c}+\Big(1+\frac{\bar{\beta}_{i-1}}{\beta_{i-1}}\Big)\frac{2}{v_f}+\frac{\xi_{i+1}}{v_fw_c}+4\Big)
    \end{align*}
    
    \item $i\in\Omega_I$
    \begin{align*}
        \frac{1}{l}\Big(\frac{2+\sqrt{2}}{w_c}+\frac{2}{v_f}+\frac{\xi_i}{v_fw_c}+5\Big)
    \end{align*}
    
    \item $i\in\Omega_O$
    \begin{align*}
        \frac{1}{l}\Big(\frac{2+\bar{\beta}_i}{w_c}+\big(2+\frac{\bar{\beta}_i}{\beta_i}+\frac{1}{\beta_i}\Big)\frac{1}{v_f}+4\Big)
    \end{align*}

    \item $i=N+\bar{i},\bar{i}\in\hat{\Omega},j\in{\Omega_I}$
    \begin{align*}
        \frac{1}{l}\Big(\frac{1}{w_c}+\frac{\xi_j}{v_fw_c}+3\Big)
    \end{align*}

    \item $i=N+N_I+\bar{i},\bar{i}\in\check{\Omega},j\in{\Omega_O}$
    \begin{align*}
        \frac{1}{l}\Big(\frac{1+\bar{\beta}_{j}}{w_c}+\Big(1+\frac{\bar{\beta}_j}{\beta_j}+\frac{1}{\beta_j}\Big)\frac{1}{v_f}+4\Big)
    \end{align*} 
\end{itemize}
So, for every $f_i,i\in[1,n]$, there exists a $\gamma_i$  such that$\mid f_i-\hat{f}_i\mid \le \gamma_i\lvert\lvert \m x - \hat{\m x}\rvert\rvert_2,\gamma_i\ge0$. Since we know that
\begin{align*}
    \lvert\lvert \m f(x,u) - \m f(\hat{x},u)\rvert\rvert_2^2=\sum_{i=1}^{n}\mid f_i-\hat{f}_i\mid^2\le\sum_{i=1}^{n}\gamma_i^2\lvert\lvert \m x - \hat{\m x}\rvert\rvert_2^2
\end{align*}
then $\gamma$ for $\m f(x,u)$ can be defined using function $h(\cdot)$ which is given by \eqref{eq:Lipschitz_constant}.
\begin{align}
    \label{eq:Lipschitz_constant}
    h\hspace{0.5mm}&(v_f,\hspace{0.5mm}\omega_c,\hspace{0.5mm}l,\hspace{0.5mm}\xi_1,\hspace{0.5mm}\xi_2,\hspace{0.5mm}\dots,\hspace{0.5mm}\xi_{N_I},\hspace{0.5mm}\beta_1,\hspace{0.5mm}\beta_2,\hspace{0.5mm}\dots,\hspace{0.5mm}\beta_{N_O})\nonumber
    \\&=\frac{1}{l}\Bigg(\Big(\frac{1+\sqrt{2}}{w_c} + \frac{1}{v_f} + \frac{\xi_2}{v_fw_c}+4\Big)^2\nonumber\\
    &+\Big(\frac{2\bar{\beta}_{N-1}}{w_c} + \Big(1+2\frac{\bar{\beta}_{N-1}}{\beta_{N-1}}\Big)\frac{1}{v_f}+4\Big)^2\nonumber\\
    &+\sum_{\substack{i\in\Omega\setminus\Omega_I\cup\Omega_O,\\i\neq\{1,N\}}}\Big(\frac{1+2\bar{\beta}_{i-1}+\sqrt{2}}{w_c}+\Big(1+\frac{\bar{\beta}_{i-1}}{\beta_{i-1}}\Big)\frac{2}{v_f}\nonumber\\
    &+\frac{\xi_{i+1}}{v_fw_c}+4\Big)^2+\sum_{i\in\Omega_I}\Big(\frac{2+\sqrt{2}}{w_c}+\frac{2}{v_f}+\frac{\xi_i}{v_fw_c}+5\Big)^2\nonumber\\
    &+\sum_{j\in\Omega_I}\Big(\frac{1}{w_c}+\frac{\xi_j}{v_fw_c}+3\Big)^2\nonumber \\
    &+\sum_{i\in\Omega_O}\Big(\frac{2+\bar{\beta}_i}{w_c}+\big(2+\frac{\bar{\beta}_i}{\beta_i}+\frac{1}{\beta_i}\Big)\frac{1}{v_f}+4\Big)^2\nonumber\\
    &+\sum_{j\in\Omega_O}\Big(\frac{1+\bar{\beta}_{j}}{w_c}+\Big(1+\frac{\bar{\beta}_j}{\beta_j}+\frac{1}{\beta_j}\Big)\frac{1}{v_f}+4\Big)^2\Bigg)^{1/2}
\end{align}

\section{Proof of Theorem \ref{thm:l_inf_theorem}}\label{apdx:obs_proof}

Let $\m z = \m Z\m e$ be the performance output for estimation error $\m e$ and $\m w$ be an unknown bounded disturbance. Let $V(\m e) = \m e^{\top}\mP \m e$ be a Lyapunov function candidate where $\mP \succ 0$. It can be shown that---for example, see \cite{Chakrabarty2016}---the estimation error dynamics \eqref{eq:est_error_dynamics} is $\mathcal{L}_{\infty}$ stable with performance level $\mu = \sqrt{\mu_0\mu_1+\mu_2}$ if there exist 
constants $\mu_0,\mu_1,\mu_2\in \mathbb{R}_{+}$ such that 
\begin{subequations}\label{eq:l_inf_lemma}
	\vspace{-0.0cm}
	\begin{align}
		%\dot{V}(\m e) &\leq -\alpha\left(V(\m e)-\mu_0 \norm{\m w}_2^2\right)\\
		\mu_0 \norm{\m w}_2^2 &\leq V(\m e) \;\Rightarrow \Delta{V}(\m e) \leq 0 \label{eq:l_inf_lemma_1}\\
		\norm{\m z}_2^2 &\leq \mu_1V(\m e)+\mu_2\norm{\m w}_2^2,\label{eq:l_inf_lemma_2}
	\end{align}
\end{subequations}
for all $t\geq 0$ where $\Delta{V}(\m e) \triangleq {V}(\m e[k+1]) - {V}(\m e[k])$.
%		 then the estimation error dynamics  \eqref{eq:est_error_dynamics} is $\mathcal{L}_{\infty}$ stable with performance level $\mu = \sqrt{\mu_0\mu_1+\mu_2}$.
First, realize that \eqref{eq:l_inf_lemma_1} holds if there exists $\alpha > 0$ such that 
\begin{align}
	\hspace{-0.2cm}{V}(\m e[k+1]) - {V}(\m e[k]) \;& \nonumber \\
	+\alpha\left(V(\m e[k])-\mu_0 \norm{\m w[k]}_2^2\right)&\leq 0 \nonumber \\
	\Leftrightarrow {\m e}^{\top}[k+1]\mP\m e[k+1] - (1-\alpha){\m e}^{\top}[k]\mP{\m e}[k]\;& \nonumber \\
	-\alpha\mu_0\m w^{\top}[k]\m w[k] &\leq 0 \nonumber 
\end{align}
which is equivalent to 
\begin{align}
	%			\bmat{\m e\\ \Delta\m f\\ \m w}^{\top} \bmat{\m \Omega &*&*\\
	%				\mP&\mO&*\\\m {B_{\mathrm{w}}}^{\top}\mP-\m {D_{\mathrm{w}}}^{\top}\mL^{\top}\mP & \mO &-\alpha\mu_0\mI} \bmat{\m e\\ \Delta\m f\\ \m w} &\leq 0\nonumber \\
	%	\Leftrightarrow 
	\bmat{\m e\\ \Delta\m f\\ \m w}^{\top}\left(\m \Pi	
	+\m\Psi^{\top}\mP^{-1}\m\Psi \right) \bmat{\m e\\ \Delta\m f\\ \m w} &\leq 0, \label{eq:l_inf_theorem_proof_1} 
\end{align}
where $\m \Pi := \mathrm{Blkdiag}\left(\left[(\alpha-1)\mP,\m O,-\alpha\mu_0\mI\right]\right)$ and $\m\Psi$ is defined as
\begin{align*}
	\m\Psi \triangleq \bmat{\mP\mA-\mP\mL\mC&\m P \m G&\mP\m {B_{\mathrm{w}}}-\mP\mL\m {D_{\mathrm{w}}}}.
\end{align*}
Since the nonlinear function $\m f(\cdot)$ is locally Lipschitz in $\m \Omega$, then it holds that
\begin{align}
	\hspace{-0.25cm}\norm{\Delta\m f[k]}_2^2 &\leq \gamma_l^2 \norm{\m e[k]}_2^2 \nonumber \\
	\Leftrightarrow 	\Delta\m f^{\top}\Delta\m f-\gamma_l^2\m e^{\top} \m e &\leq 0 \nonumber \\
	\Leftrightarrow \bmat{\m e\\ \Delta\m f\\ \m w}^{\top} \bmat{-\gamma_l^2\mI &*&*\\
		\mO&\mI&*\\\mO&\mO&\mO} \bmat{\m e\\ \Delta\m f\\ \m w} &\leq 0.\label{eq:l_inf_theorem_proof_2} 
	%\\
	%	\Leftrightarrow \mathrm{Blkdiag}\left(\left[-\gamma_l^2\mI,\mI,\mO\right]\right) &\preceq 0.\label{eq:l_inf_theorem_proof_2}
\end{align}
By applying the S-procedure to \eqref{eq:l_inf_theorem_proof_1} from \eqref{eq:l_inf_theorem_proof_2} for $\epsilon \geq 0$, we get
\begin{align}
	&\mathrm{Blkdiag}\left(\left[(\alpha-1)\mP + \epsilon\gamma_l^2\m I,-\epsilon\m I,-\alpha\mu_0\mI\right]\right) +\m\Psi^{\top}\mP^{-1}\m\Psi\preceq 0, \nonumber 
\end{align}
which, thanks to the Schur Complement, is equivalent to
\begin{align}%\label{eq:l_inf_theorem_proof_3}
	&
	\bmat{
		\renewcommand\arraystretch{1.2}
		\begin{array}{c|c}
			\mathrm{Blkdiag}\left(\left[(\alpha-1)\mP + \epsilon\gamma_l^2\m I,-\epsilon\m I,-\alpha\mu_0\mI\right]\right) \vspace*{0.0cm} & * \\ %\m\Psi^\top \\
			\hline
			\m\Psi  & -\m P
		\end{array}
	} \hspace{-0.05cm}\preceq 0.\nonumber
\end{align}
Notice that the above is equivalent to \eqref{eq:l_inf_theorem2_1}
provided that $\mY \triangleq \mP\mL$. Next, substituting $\m z = \mZ \m e$ to \eqref{eq:l_inf_lemma_2} yields
\begin{align}
	\norm{\mZ \m e}_2^2 -\mu_1V(\m e)-\mu_2\norm{\m w}_2^2&\leq 0 \nonumber\\
	\Leftrightarrow\m e^{\top}\mZ^{\top}\mZ \m e-\mu_1\m e^{\top}\mP\m e -\mu_2 \m w^{\top} \m w &\leq 0 \nonumber \\
	\Leftrightarrow \bmat{\m e\\ \m w}^{\top}\bmat{\mZ^{\top}\mZ -\mu_1\mP & \mO \\ \mO &-\mu_2\mI}\bmat{\m e\\ \m w} &\leq 0 \nonumber  \\
	\Leftrightarrow\bmat{\mZ^{\top}\mZ -\mu_1\mP & \mO \\ \mO &-\mu_2\mI} &\preceq 0 \nonumber %\\
	%	\Leftrightarrow \mathrm{Blkdiag}\left(\left[-\mu_1\m P + \m Z^\top \m Z,-\mu_2\m I\right]\right) &\preceq 0. \nonumber
\end{align}
By employing congruence transformation (provided that $\mu_1 > 0$) and applying  the Schur Complement, it is not difficult to show that the above is equivalent to \eqref{eq:l_inf_theorem2_2}. Finally, the objective function \eqref{eq:l_inf_theorem2_0} tries to minimize the performance index, which dictates how much the disturbance $\m w$ will impact the performance vector $\m z$. This implies that if optimization problem \eqref{eq:l_inf_theorem2} is solved, then the estimation error dynamics given in \eqref{eq:est_error_dynamics} is $\mathcal{L}_{\infty}$ stable with performance level $\mu = \sqrt{\mu_0\mu_1+\mu_2}$ and observer gain given as $\mL = \mP^{-1}\mY$, thus completing the proof.
\newqed
%	\end{proof}

\vspace{-0.1cm}

\end{document}